\documentclass[11pt,reqno]{amsart}

\usepackage[english]{babel}
\usepackage[utf8]{inputenc}
\usepackage[T1]{fontenc}
\usepackage{lmodern} \normalfont %to load T1lmr.fd 
\DeclareFontShape{T1}{lmr}{bx}{sc} { <-> ssub * cmr/bx/sc }{}
\usepackage{microtype}	% for improved spacing between words and letters

\usepackage{amssymb}
\usepackage{amsmath}
\usepackage{amsthm}
\usepackage{mathtools}
\usepackage{mathrsfs}
\usepackage{accents} % correct spacing for accents in sub- and superscript
\usepackage{etoolbox}
\usepackage{siunitx}
\usepackage{dsfont} % for special 1, \amthds{1}
\usepackage{graphicx}
\usepackage{color}
\usepackage[dvipsnames]{xcolor}
\usepackage{tikz}
\usetikzlibrary{calc,positioning,shapes}
\usetikzlibrary{patterns,decorations.pathmorphing,decorations.markings}
\usepackage{pgfplots}
\pgfplotsset{compat=newest}
\usepackage[margin=10pt,font=small,labelfont=bf,labelsep=endash]{caption}
\usepackage{subcaption}

\usepackage{nicematrix} %for crossing double lines in tabular

\usepackage{comment}

\usepackage{booktabs}
\usepackage{paralist}
\usepackage{soul}

\numberwithin{equation}{section}

\usepackage{algorithm}
\usepackage{algorithmicx}
\usepackage{algpseudocode}
\usepackage{cite}
\usepackage{multirow} % for multirows in tabular environment
\usepackage{enumitem} % for labelling items in enumerate
\setlist[enumerate]{label=(\roman*)}

\usepackage[colorlinks,citecolor=black,urlcolor=black]{hyperref}
\hypersetup{allcolors=black}
\usepackage[nameinlink,noabbrev]{cleveref}
\crefname{subsection}{subsection}{subsections}

\theoremstyle{plain}
\newtheorem{theorem}{Theorem}[section]

\newtheorem{remark}[theorem]{Remark}

\DeclareMathOperator*{\argmin}{arg\,min} %argmax
\newcommand{\matlab}{MATLAB\textsuperscript{\textregistered}}

\definecolor{mycolor1}{rgb}{0.00000,0.44700,0.74100}% blue
\definecolor{mycolor2}{rgb}{0.85000,0.32500,0.09800}% red
\definecolor{mycolor3}{rgb}{0.92900,0.69400,0.12500}% orange/yellow
\definecolor{mycolor4}{rgb}{0.46600,0.67400,0.18800}% green
\definecolor{mycolor5}{rgb}{0.49400,0.18400,0.55600}% purple

\title[Numerical realization of the Mortensen observer]{Numerical realization of the Mortensen observer via a Hessian-augmented polynomial approximation of the value function} %Working title
\author{Tobias Breiten${}^\dagger$ \and Karl Kunisch${}^\star$ \and Jesper Schr\"oder${}^\dagger$}
\address{${}^{\dagger}$  Institute of Mathematics MA\,{}4-4, Technical University Berlin, Stra{\ss}e des 17.~Juni 136, 10623 Berlin, Germany}
\address{${}^{\star}$ Institute of Mathematics and Scientific Computing, University of Graz, Heinrichstra\ss e 36, A-8010 Graz, Austria, and Johann Radon Institute, Austrian Academy of Sciences, A-4040 Linz, Austria}
\email{tobias.breiten@tu-berlin.de}
\email{karl.kunisch@uni-graz.at}
\email{j.schroeder@tu-berlin.de}
\date{\today}
\keywords{}

\begin{document}
	
	\begin{abstract}
		
		Two related numerical schemes for the realization of the \textit{Mortensen observer} or \textit{minimum energy estimator} for the state reconstruction of non-linear dynamical systems subject to deterministic disturbances are proposed and compared. Both approaches rely on a polynomial approximation of the value function associated with the energy of the disturbances of the system. Such an approximation is obtained via interpolation considering not only the values but also first and second order derivatives of the value function in a set of sampling points. The scheme is applied to four examples and the results are compared with the well known \textit{extended Kalman filter}.
		
	\end{abstract}
	
	\maketitle
	{\footnotesize \textsc{Keywords:} non-linear observer design, minimum energy estimation, Hamilton-Jacobi-Bellman equation, value function approximation}
	
	{\footnotesize \textsc{AMS subject classification:}} 93B53, 49M05, 49L12
	
	%%%%%%%%%%%%%%%%%%%%%%%%%%%%%%%%%%%%%%%%%%%%%%%%%%%%%%%%%%%%%%
	\section{Introduction}\label{sec: intro}
	%%%%%%%%%%%%%%%%%%%%%%%%%%%%%%%%%%%%%%%%%%%%%%%%%%%%%%%%%%%%%%
	
	We consider a non-linear dynamical system subject to linear disturbances of the form
	\begin{align*}
		\dot{x}(t) &= f(x(t)) + F v(t),~ t \in (0,T],\\
		x(0) &= x_0 + \eta,
	\end{align*}
	where $f \colon \mathbb{R}^n \rightarrow \mathbb{R}^n$, $F \in \mathbb{R}^{n,m}$, $x_0,\eta \in \mathbb{R}^n$. Here, it is assumed that $v$ and $\eta$ are unknown perturbations resulting from modeling and measurement errors or incomplete system knowledge. Our interest is to use a priori information about the dynamics $f,F,x_0$ together with (perturbed) output data $y$ on $[0,T]$ to construct a state estimate $\widehat{x}$ such that $\widehat{x}(t)\approx x(t)$ for all $t \in [0,T]$. For this purpose, we assume linear measurements subject to linear disturbances 
	\begin{align*}
		y(t) = C x(t) + \mu(t),
	\end{align*}
	where $C \in \mathbb{R}^{r,n}$ is known and $\mu \in L^2(0,T;\mathbb{R}^r)$ is unknown. The problem of reconstructing and predicting signals from partially known data has a far reaching history dating back at least to the seminal works by Wiener \cite{Wie49}, Kalman \cite{Kal60,KalB61} and Stratanovich \cite{Str60}. While the viewpoint in these articles relies on stochastic disturbances and therefore focuses on filtering theory, in this manuscript, we consider a deterministic, control theoretic perspective and aim at approximations $\widehat{x}$ characterized by an observer of the form
	\begin{align*}
		\dot{\widehat{x}}(t) &= f(\widehat{x}(t)) 
		+ L(t,\widehat{x}(t))  (y(t) - C \widehat{x}(t)),~ t \in (0,T],\\
		\widehat{x}(0) &= x_0,
	\end{align*}
	where $L$ is an appropriately chosen (non-linear) observer gain, sometimes also called output injection operator. While this observer structually resembles \textit{Luenberger} observers \cite{Lue71}, this work is concerned with a different approach. In particular, we will discuss numerical approximation schemes for the so-called \textit{Mortensen observer} or \textit{minimum energy estimator}, see \cite{Mor68,Hij80} which selects the observer gain as $L(t,\widehat{x}(t))=\nabla_{\xi\xi}^2 \mathcal{V}(t,\widehat{x}(t))^{-1}C^\top$ where $\mathcal{V}$ denotes a minimal value function for a time and output dependent Hamilton-Jacobi-Bellman (HJB) equation. It is well-known that this specific choice is a non-linear extension of and coincides with the Kalman-Bucy filter in the linear quadratic case, see \cite{Fle97}. Moreover, it can be interpreted as an \textit{optimal} estimate in the sense that it computes trajectories that are associated with perturbations of minimal energy. Though the general approach and some of its theoretical properties have been studied extensively in the literature, see, e.g., \cite{Kre79,Kre03,Moi18,BreS23,Moi23}, a numerical realization that goes beyond the case $n=3$ does not seem to be available. The main reason being that the observer gain depends on the (Hessian of the) solution to a high-dimensional HJB equation which is known to suffer from \textit{the curse of dimensionality} \cite{Bel57}. In \cite{BreKu21}, a neural network based approach was numerically investigated for the Mortensen observer for a class of non-linear oscillators corresponding to system dimensions $n=2$.
	
	At the same time, in the context of optimal feedback control where similar numerical obstacles have to be overcome, research has recently seen a tremendous progress. Without being an exhaustive list, let us for example mention value function approximations based on sparse grids \cite{GarK17,KanW17}, polynomials \cite{BreKP18,KalK18}, tensor techniques \cite{DolKK21,OstSS22} or machine learning \cite{KunW21,NakGK21,Onketal23}, respectively. Let us also mention \cite{KunW23} which discusses learning optimal feedback laws in a more general context. Here, we follow the recent work \cite{Azmetal21,KunV23} where a (sparse) polynomial learning technique has been proposed which exploits the connection between the open-loop (Pontryagin maximum principle) and the closed-loop (HJB equation) point of view by first sampling the value function and its gradient and subsequently solving a (sparsity promoting) least squares problem. Let us emphasize that, despite its conceptual similarity to those works, value function approximations tailored to the Mortensen observer differ considerably from the optimal feedback framework for several reasons. First of all, the HJB equation does explicitly depend on the output trajectory $y$ resulting in a generic time-dependent source term which has to be taken into account. In particular, we will point out that already in the linear quadratic case, the structure of the value function with respect to the temporal and the spatial components differ significantly, complicating the choice of appropriate polynomial basis functions and their degrees. 
	Furthermore, in contrast to an optimal feedback law, the computation of the observer gain involves the inverse of the Hessian which is prone to numerical errors if, e.g., finite difference approximations were used. Finally, since the aforementioned works on feedback control deal with the problem of optimal stabilization (or tracking to zero), sampling of the value function is naturally done within a sufficiently large region around the origin. On the contrary, for the minimum estimation problem, we will see that an accurate approximation is particularly important along the observer trajectory which of course is unknown a priori. 
	
	With this in mind, our contributions are as follows:
	\begin{itemize}
		\item Based on the theoretical characterization of the Mortensen observer by means of an optimal control problem, we derive a \textit{differential Riccati equation} (DRE) (depending on the adjoint state) which provides Hessian information of the value function. In particular, this allows us to appropriately enrich the data set for our least squares problem.
		\item By combining quasi-random Halton sampling around the trajectory of an extended Kalman filter, we propose a numerically efficient data generation procedure. This strategy is theoretically motivated by the fact that the extended Kalman filter can be understood as a first-order approximation to the Mortensen observer.
		\item Resorting to its original formulation in \cite{Mor68}, in addition to the Luenberger-type formulation of the Mortensen observer, we also address an alternative formulation which is based on the pointwise (in time) minimization of the value function over the state space and numerically investigate its performance when compared to the first formulation. In a similar spirit, we also compare the extended Kalman filter to our approximations of the Mortensen observer with particular emphasis on systems whose value functions exhibit a spatially strongly non-quadratic behavior.
		\item With regard to its applicability in the context of larger state space dimensions, we utilize a hyperbolic cross approximation which allows us to compute HJB-based observers up to dimensions $n=40$. 
	\end{itemize}
	
	\vspace{5mm}
	
	\textbf{Notation.}
	If not mentioned otherwise, $\Vert \cdot \Vert$ will denote the Euclidean norm on $\mathbb{R}^d$, where the dimension $d$ varies. The associated scalar product is denoted by $\langle \cdot, \cdot \rangle$. Further we denote by $I_d$ the identity matrix of dimension $d$. The space of matrices with $d_1$ rows and $d_2$ columns and real entries will be denoted by $\mathbb{R}^{d_1,d_2}$. If not mentioned otherwise it is equipped with the Frobenius norm. For $1 \leq p \leq \infty $ we denote by $L^p(0,T;\mathbb{R}^d)$ the Lebesgue spaces, while $H^1(0,T;\mathbb{R}^d)$ denotes the Sobolev space of functions with a first weak derivative in $L^2(0,T;\mathbb{R}^d)$. All mentioned function spaces are equipped with their respective standard norms. For two Banach spaces $X$ and $Y$ we denote by $L(X,Y)$ the space of linear and bounded mappings from $X$ to $Y$. For a function $f \colon X \to Y$ its Fr\'{e}chet derivative is denoted by $Df$. In case $X$ and $Y$ are finite dimensional, we identify $Df$ with the Jacobian matrix. 
	
	%%%%%%%%%%%%%%%%%%%%%%%%%%%%%%%%%%%%%%%%%%%%%%%%%%%%%%%%%%%%%%
	\section{The Mortensen observer}\label{sec: MorObs}
	%%%%%%%%%%%%%%%%%%%%%%%%%%%%%%%%%%%%%%%%%%%%%%%%%%%%%%%%%%%%%%
	
	In this section we briefly recall the motivation and definition of the well-known \textit{Mortensen observer} and the essential concepts needed for its construction. A more thorough discussion can  be found in \cite{BreKu21}, for example, see also the original work \cite{Mor68}. The strategy of the approach at hand is to minimize the energy of the disturbances in the system. The mathematical derivation is based on a specific optimal control problem. For fixed $t \in (0,T]$, $\xi \in \mathbb{R}^n$, and $y \in L^2(0,T;\mathbb{R}^r)$ it reads
	\begin{equation}\label{eq: OCP}
		\begin{aligned}
			\inf_{\substack{x \in H^1(0,t;\mathbb{R}^n) \\ v \in L^2(0,t;\mathbb{R}^m)}}
			J(x,v;t) &\coloneqq
			\frac{1}{2} \Vert x(0) - x_0 \Vert^2
			+ \frac{1}{2} \int_0^t \Vert v(s) \Vert^2 + \alpha \Vert y(s) - C x(s) \Vert^2 \,\mathrm{d}s,\\
			\text{s.t. } e(x,v;t,\xi) &\coloneqq (\dot{x} - f(x) - Fv, x(t) - \xi) = 0.
		\end{aligned}
	\end{equation}
	For $\xi \in \mathbb{R}^n$ the corresponding minimal value function is defined by
	\begin{equation}\label{eq: ValFun}
		\mathcal{V}(t,\xi) \coloneqq 
		\inf_{\substack{(x,v) \text{ s.t.} \\ e(x,v;t,\xi) = 0}}J(x,v;t),~ t \in (0,T],~~~~~~~~
		\mathcal{V}(0,\xi) \coloneqq 
		\frac{1}{2} \Vert \xi - x_0 \Vert^2.
	\end{equation}
	From the point of view of optimal control theory problem \eqref{eq: OCP} differs from the classic literature with respect to the boundary condition which is given at the right boundary of the time interval instead of the left one. Indeed the system evolves backwards in time, which has a crucial effect on the numerical treatment. Systems that are stable when considered forward in time might turn unstable when considered as evolving backwards in time. This issue will be discussed in more detail in our examples.
	
	We recall that $\mathcal{V}$ can be characterized as the solution of the \textit{Hamilton-Jacobi-Bellman equation} associated with \eqref{eq: OCP}. Since the state equation evolves backwards in time, in our case the HJB equation is given as a forward equation and reads
	\begin{equation}\label{eq: HJB}
		\begin{aligned}
			\partial_t \mathcal{V}(t,\xi) &=
			-\nabla_\xi \mathcal{V}(t,\xi)^\top f(\xi) 
			- \frac{1}{2} \Vert F^\top \nabla_\xi \mathcal{V}(t,\xi) \Vert^2
			+ \frac{\alpha}{2} \Vert y(t) - C\xi \Vert^2,\\
			\mathcal{V}(0,\xi) &= \frac{1}{2} \Vert \xi - x_0 \Vert^2.
		\end{aligned}
	\end{equation}
	Depending on the regularity of the underlying control problem one might have to consider \textit{viscosity solutions} instead of classical ones, c.f.~\cite{FleS06}. In \cite{BreS23} conditions were given which ensure that $\mathcal{V}$ is in fact space-time $C^1$-regular in a neighborhood of the model.
	
	The Mortensen observer is defined via a pointwise minimization of the value function, i.e.,
	\begin{equation}\label{eq: observerDef}
		\widehat{x}_{\mathrm{M}}(t) \coloneqq \argmin\limits_{\xi \in \mathbb{R}^n} \mathcal{V}(t,\xi).
	\end{equation}
	Under appropriate assumptions it can be characterized as a solution of the observer equation 
	\begin{equation}\label{eq: observerEq}
		\begin{aligned}
			\dot{\widehat{x}}_{\mathrm{M}}(t) &= f(\widehat{x}_{\mathrm{M}}(t)) 
			+ \alpha \nabla_{\xi\xi}^2 \mathcal{V}(t,\widehat{x}_{\mathrm{M}}(t))^{-1} C^\top (y(t) - C \widehat{x}_{\mathrm{M}}(t)),~ t \in (0,T],\\
			\widehat{x}_{\mathrm{M}}(0) &= x_0,
		\end{aligned}
	\end{equation}
	c.f.~\cite{BreS23}.
	
	In the derivation and motivation of our numerical strategies we assume that the underlying system and functions fulfill all necessary assumptions needed to ensure well-definedness of \eqref{eq: observerDef} and well-posedness of \eqref{eq: observerEq}. This specifically includes sufficient smoothness of the value function. A theoretical result on these issues for systems with a quadratic right-hand side can be found in \cite{BreS23}, where existence of a solution $\widehat{x}_{\mathrm{M}} \in H^1(0,T;\mathbb{R}^n)$ was established.
	
	\vspace{5mm}
	
	\textbf{Relation to the Kalman filter.}
	In the case of linear dynamics, i.e., $f(\xi) = A \xi$ for some $A \in \mathbb{R}^{n,n}$, it can be shown that the value function is quadratic in $\xi$ and explicitly given by
	\begin{equation}\label{eq: KalVF}
		\mathcal{V}(t,\xi) = 
		\frac{1}{2} (\xi - \widehat{x}(t))^\top \Sigma(t)^{-1} (\xi - \widehat{x}(t))
		+ \frac{\alpha}{2} \int_0^t \Vert y(s) - C \widehat{x}(s) \Vert^2 \,\mathrm{d} s,
	\end{equation}
	where 
	\begin{equation}\label{eq: Kalman}
		\begin{aligned}
			\dot{\Sigma}(t) &= A \, \Sigma(t) + \Sigma(t) \, A^\top - \alpha \Sigma(t) \, C^\top C \, \Sigma(t) + FF^\top,~ \Sigma(0) = I_n,\\
			\dot{\widehat{x}}(t) &= A \widehat{x}(t) 
			+ \alpha \Sigma(t) \, C^\top (y(t) - C \widehat{x}(t)),~
			\widehat{x}(0) = x_0,
		\end{aligned}
	\end{equation}
	c.f.~\cite{BreKu21}. Note that $\Sigma(t)$ coincides with the inverted Hessian $\nabla_{\xi\xi}^2 \mathcal{V}(t,\widehat{x}(t))^{-1}$ and is given as the solution of a differential Riccati equation, while $\widehat{x}$ is given as the solution of the observer equation which agrees with the formulation in \eqref{eq: observerEq}. System \eqref{eq: Kalman} is the deterministic version of the widely used \textit{Kalman-Bucy filter}\cite{Kal60,KalB61}.

	\vspace{5mm}
	
	\textbf{The extended Kalman filter.}
	In engineering practice observers for non-linear problems are frequently  based on the use of the Kalman filter on a linearization of the system at the observer trajectory. This leads to the so-called extended Kalman filter which for our setting is formulated as
	
	\begin{equation}\label{eq: ext Kalman}
		\begin{aligned}
			\dot{\Sigma}(t) &= \mathrm{D}f(\widehat{x}_{\mathrm{K}}(t)) \, \Sigma(t) + \Sigma(t) \, \mathrm{D}f(\widehat{x}_{\mathrm{K}}(t))^\top 
			- \alpha \Sigma(t) \, C^\top C \, \Sigma(t) + FF^\top,~ \Sigma(0) = I_n, \\
			\dot{\widehat{x}}_{\mathrm{K}}(t) &= f(\widehat{x}_{\mathrm{K}}(t)) + \alpha \Sigma(t) \, C^\top ( y(t) - C \widehat{x}_{\mathrm{K}}(t) ),~ \widehat{x}_{\mathrm{K}}(t) = x_0.
		\end{aligned}
	\end{equation}
	
	In contrast to \eqref{eq: Kalman} here the differential Riccati equation is coupled with the observer equation. For our non-linear examples we will employ an implementation of this extended Kalman filter as a first approximation of the Mortensen observer. While the extended Kalman filter is a powerful technique, we shall also present an example demonstrating that the state reconstruction $\widehat{x}_\mathrm{M}$ based on Mortensen can significantly differ from the state reconstruction $\widehat{x}_\mathrm{K}$ obtained by means of the extended Kalman filter.
	
	We further point out that the extended Kalman filter and the Mortensen observer are closely related. Both are based on feeding the observation defect back into the model using gains that are characterized via respective differential Riccati(-type) equations. To see that the latter holds for the Mortensen observer consider the second order spatial derivative of the HJB \eqref{eq: HJB}
	\begin{equation}\label{eq: HJB 2 der}
		\begin{aligned}
			\partial_t \nabla_{\xi\xi}^2 \mathcal{V}(t,\xi)
			&= - \nabla_{\xi\xi}^2 \mathcal{V}(t,\xi) \mathrm{D}f(\xi)
			- \mathrm{D}f(\xi)^\top \nabla_{\xi\xi}^2 \mathcal{V}(t,\xi)\\
			&- \nabla_{\xi\xi}^2 \mathcal{V}(t,\xi) FF^\top \nabla_{\xi\xi}^2 \mathcal{V}(t,\xi)
			- \mathrm{D} \left( \mathrm{D}f(\xi)^\top \right) \nabla_\xi \mathcal{V}(t,\xi)
			+ \alpha C^\top C\\
			&- \nabla_{\xi^3}^3 \mathcal{V}(t,\xi) \left( f(\xi) + FF^\top \nabla_\xi \mathcal{V}(t,\xi) \right).
		\end{aligned}
	\end{equation}
	Due to \eqref{eq: observerDef} it holds $\nabla_\xi \mathcal{V}(t,\widehat{x}_{\mathrm{M}}(t)) = 0$ for all $t \in [0,T]$. With \eqref{eq: HJB 2 der} it follows
	\begin{equation}\label{eq: Riccati4Hess}
		\begin{aligned}
			\partial_t \nabla_{\xi\xi}^2 \mathcal{V}(t,\widehat{x}_{\mathrm{M}}(t))
			&= - \nabla_{\xi\xi}^2 \mathcal{V}(t,\widehat{x}_{\mathrm{M}}(t)) \mathrm{D}f(\widehat{x}_{\mathrm{M}}(t))
			- \mathrm{D}f(\widehat{x}_{\mathrm{M}}(t))^\top \nabla_{\xi\xi}^2 \mathcal{V}(t,\widehat{x}_{\mathrm{M}}(t))\\
			&- \nabla_{\xi\xi}^2 \mathcal{V}(t,\widehat{x}_{\mathrm{M}}(t)) FF^\top \nabla_{\xi\xi}^2 \mathcal{V}(t,\widehat{x}_{\mathrm{M}}(t))
			+ \alpha C^\top C \\
			&- \nabla_{\xi^3}^3 \mathcal{V}(t,\widehat{x}_{\mathrm{M}}(t)) f(\widehat{x}_{\mathrm{M}}(t)).
		\end{aligned}
	\end{equation}
	We define $\Pi(t) = \nabla_{\xi\xi}^2 \mathcal{V}(t,\widehat{x}_{\mathrm{M}}(t))^{-1}$ and multiply \eqref{eq: Riccati4Hess} with $\Pi(t)$ from the left and the right. Subsequently the Mortensen observer can be characterized via
	\begin{equation}\label{eq: MortViaRicc}
		\begin{aligned}
			\dot{\Pi}(t) &= 
			\mathrm{D}f(\widehat{x}_{\mathrm{M}}(t)) \Pi(t) + \Pi(t) \mathrm{D}f(\widehat{x}_{\mathrm{M}}(t))^\top 
			- \alpha \Pi(t) C^\top C \Pi(t) + F F^\top \\
			&+ \Pi(t) \nabla_{\xi^3}^3 \mathcal{V}(t,\widehat{x}_{\mathrm{M}}(t)) f(\widehat{x}_{\mathrm{M}}(t)) \Pi(t),~ \Pi(0) = I_n,\\
			\dot{\widehat{x}}_{\mathrm{M}}(t) &= f(\widehat{x}_{\mathrm{M}}(t)) + \alpha \Pi(t) \, C^\top ( y(t) - C \widehat{x}_{\mathrm{M}}(t) ),~ \widehat{x}_{\mathrm{M}}(t) = x_0,
		\end{aligned}
	\end{equation}
	and thus a comparison with \eqref{eq: ext Kalman} shows that the gains for the extended Kalman filter and the Mortensen observer are characterized via Riccati(-type) equations that differ only by the summand involving the third derivative of the value function. 
	
	However, since usually the third derivative $\nabla_{\xi^3}^3 \mathcal{V}$ is not available in practical computations, this characterization does not offer any advantages when approximating the Mortensen observer.
	
	%%%%%%%%%%%%%%%%%%%%%%%%%%%%%%%%%%%%%%%%%%%%%%%%%%%%%%%%%%%%%%
	\section{Polynomial approximation of the value function}\label{sec: AppVF}
	%%%%%%%%%%%%%%%%%%%%%%%%%%%%%%%%%%%%%%%%%%%%%%%%%%%%%%%%%%%%%%
	
	In this work we propose a scheme for the numerical approximation of the Mortensen observer for non-linear systems. We will consider two different strategies, namely construction by pointwise minimization of the value function according to \eqref{eq: observerDef} and by solving the observer equation \eqref{eq: observerEq}. For both approaches an appropriate approximation of the value function $\mathcal{V}$ is essential.
	
	Due to the \textit{curse of dimensionality} solving the HJB equation \eqref{eq: HJB} for $\mathcal{V}$ is extremely challenging for systems of medium and high state dimension. Instead, we sample the value function and use linear regression to obtain a polynomial approximation of $\mathcal{V}$. This approach is  inspired by \cite{Azmetal21}, where the samples of the value function are augmented by samples of its gradient. We extend this idea to also sampling the Hessian matrices. Indeed, the structure of the observer equation \eqref{eq: observerEq} calls for an accurate approximation of the Hessian of $\mathcal{V}$.  Here the Hessian samples turn out to be very helpful. 
	
	In the following we describe the steps necessary for the approximation of the value function. First we discuss the generation of a data set for training by solving open-loop optimal control problems and differential Riccati equations. The former yields the values and first order derivatives of the value function, while the latter provides the Hessian matrices. Subsequently we shall present the construction of an appropriate polynomial basis using Chebyshev polynomials and the hyperbolic cross index set. Approximation properties of hyperbolic cross based polynomials are well analyzed in the literature. We refer to, e.g., \cite[Section 4.2]{DTU18} where convergence rates in terms of Sobolev- and Besov space norms are provided. Finally we are prepared to construct a polynomial approximation of the value function $\mathcal{V}$ using linear regression.
	
	%%%%%%%%%%%%%%%%%%%%%%%%%%%%%%%%%%%%%%%%%%%%%%%%%%%%%%%%
	\subsection{Generating a dataset}\label{subsec: data}~\\
	%%%%%%%%%%%%%%%%%%%%%%%%%%%%%%%%%%%%%%%%%%%%%%%%%%%%%%%%
	
	First we generate a data set of the form 
	\begin{equation}\label{eq: sampl data}
		\left\{ (t^i,\xi^i), \mathcal{V}(t^i,\xi^i), \nabla_\xi \mathcal{V}(t^i,\xi^i), \nabla_{\xi\xi}^2 \mathcal{V}(t^i,\xi^i) \right\}_{i=1}^N.
	\end{equation}
	Our strategy of choosing the sample points $(t^i,\xi^i)$ is presented below. First we discuss how the value function and its first and second order derivatives are evaluated for such a given point. For this purpose let $(t^*,\xi^*) \in (0,T] \times \mathbb{R}^n$ denote a generic sampling point. We numerically solve the associated open-loop optimal control problem \eqref{eq: OCP} using a gradient descent approach based on \textit{Pontryagin's maximum principle}. According to this principle an optimal triple $(\Bar{x},\Bar{v},p) \in H^1(0,t^*;\mathbb{R}^n) \times L^2(0,t^*;\mathbb{R}^m) \times H^1(0,t^*;\mathbb{R}^n)$ consisting of a state trajectory, control and adjoint state satisfies the first order necessary optimality condition:
	\begin{align}
		\dot{\bar{x}} &= f(\bar{x}) + F\bar{v},~\bar{x}(t^*) = \xi^*, \label{algn: PMP state}\\
		- \dot{p} &= \mathrm{D}f(\bar{x})^\top p - \alpha C^\top (y - C\bar{x}),~p(0) =  x_0 - \bar{x}(0),\label{algn: PMP adj}\\
		\bar{v} + F^\top p &= 0. \label{algn: PMP OC}
	\end{align}
	Note that $\bar{v} + F^\top p$ coincides with the gradient of the reduced cost functional $\tilde{J}$ which is given by
	\begin{equation*}
		\tilde{J}(v;t) = J(x(v),v;t), 
	\end{equation*}
	where $x(v)$ denotes the solution of the state equation with control $v$, i.e., $e(x(v),v;t,\xi) = 0$ holds. 
	
	With a representation of the gradient at hand we can set up the gradient descent scheme. We denote by $v_j$ the control chosen in the $j$-th iteration, $x_j$ and $p_j$ are the corresponding solutions of the state and adjoint equations, respectively. We express the corresponding gradient as $\mathcal{G}_j = v_j + F^\top p_j$.
	
	For choosing a stepsize in each iteration of the gradient descent we opt for a combination of the \textit{Barzilai-Borwein} step-size-control \cite{BarB88,AzmK20} and a specific line search scheme. In the $k$-th iteration the initial stepsize is set to ${\sigma_k}^{-1}$, where
	\begin{equation}\label{eq: BarBor step}
		\sigma_k \coloneqq 
		\frac{\left\langle \mathcal{S}_{k-1} , \mathcal{Y}_{k-1} \right\rangle}
		{\left\langle \mathcal{S}_{k-1} , \mathcal{S}_{k-1} \right\rangle}, ~\text{ for even } k, \text{ and }~~~~~~~~~
		\sigma_k \coloneqq
		\frac{\left\langle \mathcal{Y}_{k-1} , \mathcal{Y}_{k-1} \right\rangle}
		{\left\langle \mathcal{S}_{k-1} , \mathcal{Y}_{k-1} \right\rangle}, ~\text{ for odd } k,
	\end{equation}
	and $\mathcal{S}_{k-1} \coloneqq v_k - v_{k-1} $ and $\mathcal{Y}_{k-1} \coloneqq \mathcal{G}_{k} - \mathcal{G}_{k-1}$. The stepsize $h_k$ used in the $k$-th iteration is then determined by an application of the non monotone line search introduced in \cite{DaiZ01}.
	
	%%%%%%%%%%%%%%%%%%%%%%%%%%%%%%%%%%%%%%%%%%%%
	%%%%%%%%ALGORITHM 1: OPEN LOOP SOLVER%%%%%%%
	%%%%%%%%%%%%%%%%%%%%%%%%%%%%%%%%%%%%%%%%%%%%
	
	\begin{algorithm}
		\caption{Gradient descent for open-loop optimal control problem}\label{alg: OLS}
		\begin{algorithmic}
			\Require Initial controls $v_{-1}$ and $v_0$, tolerances $\varepsilon_{\mathrm{rel}}$, $\varepsilon_{\mathrm{abs}}> 0$.
			\Ensure Optimal triple $(\bar{x},\bar{v},p)$
			\State Set $k$ = 0.
			\State Compute $x_{-1}$, $x_0$ via \eqref{algn: PMP state}.
			\State Compute $p_{-1}$, $p_0$ via \eqref{algn: PMP adj}.
			\State Set $\mathcal{G}_{-1} = v_{-1} + F^\top p_{-1}$, $\mathcal{G}_{0} = v_{0} + F^\top p_{0}$ 
			\While{ $ \Vert \mathcal{G}_k \Vert / \Vert \mathcal{G}_{-1} \Vert > \varepsilon_{\mathrm{rel}} $ and $\Vert \mathcal{G}_k \Vert > \varepsilon_{\mathrm{abs}}$ }
			\State Compute $\sigma_k$ according to \eqref{eq: BarBor step}.
			\State Obtain stepsize $h_k$ by non monotone line search according to \cite{DaiZ01} starting with $\sigma_k^{-1}$. 
			\State Set $v_{k+1} = v_k - h_k \mathcal{G}_k$.
			\State Compute $x_{k+1}$ via \eqref{algn: PMP state}.
			\State Compute $p_{k+1}$ via \eqref{algn: PMP adj}.
			\State Set $\mathcal{G}_{k+1} = v_{k+1} + F^\top p_{k+1}$.
			\State Set $k = k + 1$.
			\EndWhile
		\end{algorithmic}
	\end{algorithm}
	
	\Cref{alg: OLS} allows the computation of approximations of the optimal state, control, and adjoint state $(\bar{x},\bar{v},p)$ for a given $(t^*,\xi^*)$. The value of the cost functional can be computed as $\mathcal{V}(t^*,\xi^*) = J(\bar{x},\bar{v};t^*)$. Via a feedback rule one obtains the gradient as $\nabla_\xi \mathcal{V}(t^*,\xi^*) = -p(t^*)$ without any additional computational cost. The data point is completed by a computation of the Hessian which is achieved in the following manner.
	
	\vspace{5mm}
	
	\textbf{Evaluation of the Hessian in sample points.}
	First we show that $\Xi(s) = \nabla_{\xi\xi}^2 \mathcal{V}(s,\bar{x}(s))$ satisfies a specific differential Riccati equation. By the verification theorem the optimal control can be characterized via the feedback rule \cite{CanF91}
	\begin{equation}\label{eq: feedback}
		\bar{v}(s) = -F^\top p(s) = F^\top \nabla_\xi \mathcal{V}(t,\bar{x}(s)),~~~s\in [0,t^*].
	\end{equation}
	With the chain rule it follows that
	\begin{equation*}
		\tfrac{\mathrm{d}}{\mathrm{d}s} \Xi(s)
		= \partial_t \nabla_{\xi\xi}^2 \mathcal{V}(s,\bar{x}(s)) + \nabla_{\xi^3}^3 \mathcal{V}(s,\bar{x}(s)) \dot{\bar{x}}(s).
	\end{equation*}
	Inserting \eqref{eq: HJB 2 der} with $\xi = \bar{x}(s)$ yields
	\begin{equation*}
		\begin{aligned}
			\dot{\Xi}(s) &= 
			- \nabla_{\xi\xi}^2 \mathcal{V}(s,\bar{x}(s)) \mathrm{D}f(\bar{x}(s))
			- \mathrm{D}f(\bar{x}(s))^\top \nabla_{\xi\xi}^2 \mathcal{V}(s,\bar{x}(s))\\
			&- \nabla_{\xi\xi}^2 \mathcal{V}(s,\bar{x}(s)) FF^\top \nabla_{\xi\xi}^2 \mathcal{V}(s,\bar{x}(s))
			- \mathrm{D} \left( \mathrm{D}f(\bar{x}(s))^\top \right) \nabla_\xi \mathcal{V}(s,\bar{x}(s))
			+ \alpha C^\top C\\
			&- \nabla_{\xi^3}^3 \mathcal{V}(s,\bar{x}(s)) \left( f(\bar{x}(s)) + FF^\top \nabla_\xi \mathcal{V}(s,\bar{x}(s)) - \dot{\bar{x}}(s) \right).
		\end{aligned}
	\end{equation*}
	With \eqref{eq: feedback} it follows that the last summand on the right hand side vanishes and thus
	\begin{equation}\label{eq: DRE samp}
		\begin{aligned}
			\dot{\Xi}(s) &= 
			- \Xi(s) \mathrm{D}f(\bar{x}(s)) 
			- \mathrm{D}f(\bar{x}(s))^\top \Xi(s)
			- \Xi(s) FF^\top \Xi(s)\\
			&+ \mathrm{D} \left( \mathrm{D} f(\bar{x}(s))^\top \right)  p(s)
			+ \alpha C^\top C,\\
			\Xi(0) &= I_n.
		\end{aligned}
	\end{equation}
	The initial condition follows directly from the definition in \eqref{eq: ValFun}. After solving this DRE we obtain $\nabla_{\xi\xi}^2 \mathcal{V}(t^*,\xi^*) = \nabla_{\xi\xi}^2 \mathcal{V}(t^*,\bar{x}(t^*)) = \Xi(t^*)$ and the data point is complete.
	
	\begin{remark}
		It is possible to augment the sampled data by the time derivative. Since we already have access to the gradient, the time derivative $\partial_t \mathcal{V}$ can be computed via the HJB equation \eqref{eq: HJB} without any noteworthy additional cost. However, this did not yield any convincing advantages in our experiments and it is therefore omitted in this work.  
	\end{remark}
	
	\vspace{5mm}
	
	\textbf{Choosing sampling points.}
	For the generation of the data set \eqref{eq: sampl data} an appropriate choice of sampling points $(t_i,\xi_i)_{i=1}^N$ is essential. Note that we will only require evaluations of the value function and its derivatives close to the observer trajectory $\widehat{x}_{\mathrm{M}}$. For the approach based on solving the observer equation this can immediately be seen in \eqref{eq: observerEq}. The same holds when minimizing the value function over the space variable if one assumes to have access to a sufficiently good initial candidate for the minimization, c.f., \Cref{sec: realObs}. 
	
	The extended Kalman filter offers an approximation $\widehat{x}_{\mathrm{K}}$ of the Mortensen observer at comparatively small computational cost. Therefore we shall sample the value function locally around the observer trajectory $\widehat{x}_{\mathrm{K}}$. This is done in the following manner:
	
	First $N_{\mathrm{Time}}$ time sample points are determined as the Chebyshev nodes in $[0,T]$, i.e.,
	\begin{equation*}
		t_k = \tfrac{1}{2} \, T + \tfrac{1}{2} \, T \cos \left( \frac{2 \, (N_{\mathrm{Time}} - k + 1) - 1 }{2 \,  \pi \, N_{\mathrm{Time}} } \right),~~~~~~k=1,..., N_{\mathrm{Time}}.
	\end{equation*}
	These are precisely the roots of the Chebyshev polynomial of degree $N_{\mathrm{Time}}$ rescaled to the domain $[0,T]$. In the following every time sampling point $t_k$ will be individually completed by appropriate spatial sampling points. 
	
	Turning to the spatial sampling, for each fixed $k$ the Kalman observer trajectory $\widehat{x}_{\mathrm{K}}$ is evaluated in $t_k$. The spatial variable $\xi \in \mathbb{R}^n$ is sampled quasi-randomly from a hyperrectangle $Q^k$ around $\widehat{x}_{\mathrm{K}}(t_k)$ of the form
	\begin{equation}\label{eq: samp rec}
		Q^k 
		= \bigtimes_{i=1}^n \left[ \widehat{x}_{\mathrm{K}}(t_k)_i - r_{k,i},\widehat{x}_{\mathrm{K}}(t_k)_i + r_{k,i} \right],
	\end{equation}
	where the specific side lengths $r_{k,i} > 0$ are chosen individually for every example and depend on $\widehat{x}(t_k)$, c.f.~\Cref{sec: num}. We take $N_{\mathrm{Space}}$ spatial samples $(\xi_h^{t_k})_{h=1}^{N_{\mathrm{Space}}}$ from $Q^k$ using Halton quasi-random sequences\footnote{\url{https://de.mathworks.com/help/stats/generating-quasi-random-numbers.html}}. Finally the set of sampling points is given as
	\begin{equation*}
		(t_i,\xi_i)_{i=1}^N = \bigcup\limits_{k=1,h=1}^{N_{\mathrm{Time}},N_{\mathrm{Space}}} \left\{ (t_k,\xi_h^{t_k}) \right\}
	\end{equation*}
	consisting of $N = N_{\mathrm{Time}} \, N_{\mathrm{Space}}$ sample points.
	
	Let us note that initially we used tensorized Chebyshev nodes as spatial samples which worked out for examples of lower dimension. In higher dimensions, however, this is not feasible because the tensorization quickly results in a very large number of samples. An attempt of using random spatial samples uniformly distributed in $Q^k$ did not yield the desired results as samples would tend to cluster instead of filling out the sampling domain evenly. Eventually we decided to construct the samples in a quasi-random fashion as it is done in \cite{Azmetal21}.
	
	%%%%%%%%%%%%%%%%%%%%%%%%%%%%%%%%%%%%%%%%%%%%%%%%%%%%%%%%
	\subsection{Constructing a polynomial basis}\label{subsec: PolBas}~\\
	%%%%%%%%%%%%%%%%%%%%%%%%%%%%%%%%%%%%%%%%%%%%%%%%%%%%%%%%
	
	In this section we present the construction of the polynomial basis which we use for the approximation of the value function. The $(n+1)$-dimensional basis polynomials will be given as products of one-dimensional Chebyshev polynomials representing the time variable and each spatial variable separately.
	
	First an appropriate domain
	$$\mathcal{D} = \mathcal{D}_{\text{Time}} \times \mathcal{D}_{\text{Space}}  \subset \mathbb{R} \times \mathbb{R}^n$$
	needs to be determined. Here
	$$
	\mathcal{D}_{\text{Space}} = \bigtimes\limits_{i=1}^n \mathcal{D}_i \subset \mathbb{R}^n,
	$$
	and $\mathcal{D}_{\text{time}}$ and $\mathcal{D}_i$ are the domains corresponding to the Chebyshev polynomials in the time variable and in the $i$-th spatial variable respectively. While the choice $\mathcal{D}_{\text{time}} = [0,T]$ is clear, the choice for the spatial domains is less obvious. We aim at choosing domains as small as possible while still ensuring that for any $t \in [0,T]$ the set $\mathcal{D}_{\text{Space}}$ contains all spatial samples $(\xi_i)_{i=1}^N$ and the evaluations of the (unknown) observer trajectory $\widehat{x}_\mathrm{M}(t)$. Using the notation from \eqref{eq: samp rec} we define for all $i \in \{ 1,...,n \}$
	$$
	\mathcal{D}_i \coloneqq \left[\min_{k = 1,...,M_{\text{Time}}} \widehat{x}_\mathrm{K}(t_k)_i - r_{k,i}, \max_{k = 1,...,M_{\text{Time}}} \widehat{x}_\mathrm{K}(t_k)_i + r_{k,i} \right].
	$$
	
	By rescaling the Chebyshev polynomials (originally defined on $(-1,1)$) to the respective domains before normalizing, for the $i$-th spatial variable we obtain a one-dimensional orthonormal basis of $L^2(\mathcal{D}_i)$ and denote the basis functions by $(\phi_k^i)_{k = 0}^\infty$. Analogously we construct a basis of $L^2([0,T])$ denoted by $(\phi_k^{\text{Time}})_{k = 0}^\infty$. These functions are used to define the elements of the tensor-product basis of $L^2(\mathcal{D})$ via
	\begin{equation*}
		\Phi_i(t,\xi) \coloneqq \phi_{i_1}^{\text{Time}}(t) \, \prod_{j = 1}^n \phi_{i_{j+1}}^j (\xi_j), ~~~~~\text{ with } i = (i_1,...,i_{n+1}) \in \mathbb{N}_0^{n+1},~ \xi = (\xi_1,...,\xi_n). 
	\end{equation*}
	Now assuming $\mathcal{V} \in L^2\left([0,T] \times \mathcal{D}_{\text{space}}\right) \cap L^\infty\left([0,T] \times \mathcal{D}_{\text{space}}\right) $ the value function can be expanded into a series of the form
	\begin{equation*}
		\mathcal{V}(t,\xi) = \sum_{i \in \mathbb{N}_0} \theta_i \Phi_i(t,\xi),
	\end{equation*}
	with $\theta_i = \langle \mathcal{V} , \Phi_i \rangle_{L^2(\mathcal{D})}$. This motivates the approximation of the value function in a polynomial basis of the form $\{ \Phi_i \}_{i \in \mathcal{J}} $, where $\mathcal{J} \subset \mathbb{N}_0^{n+1}$ is an index set of finite cardinality.
	
	A proper choice of $\mathcal{J}$ is crucial to the accuracy of the approximation and to the numerical feasability of its computation. Especially for large spatial dimensions $n$ one needs to pay attention to the choice of indices to make sure that the set of basis polynomials does not grow too large. Since the dependence of the value function on the time variable and on the spatial variables are fundamentally different, we choose to treat them separately. Our treatment of the spatial polynomials is heavily inspired by what was done in \cite{Azmetal21}, see also \cite{ABB}, and the monograph \cite{DTU18}, therefore we will only briefly summarize the strategy here. Namely we make use of the hyperbolic cross index set to treat the spatial polynomials and obtain a tensor-product basis of $\mathcal{D}_{\text{Space}}$. For a fixed number $s \in \mathbb{N}$ the index set in question is defined as
	\begin{equation*}
		\mathcal{J}_{\text{Space}}(s) \coloneqq \left\{ i = (i_1,...,i_n) \in \mathbb{N}_0^n~ \colon ~ \prod_{j=1}^n (i_j + 1) \leq s + 1 \right\}.
	\end{equation*}
	The elements of this basis will be fully tensorized with the polynomials $(\phi_k^{\text{Time}})_{k=0}^{d_{\text{Time}}}$ up to some prescribed maximal degree $d_{\text{Time}} \in \mathbb{N}_0$. We obtain the polynomial basis represented by the index set
	$$
	\mathcal{J} = \mathcal{J}(d_{\text{Time}},s) \coloneqq \left\{ i = (i_{\text{Time}},i_{\text{Space}}) ~ \colon ~ i_{\text{Time}} \in \{ 0,..., d_{\text{Time}} \}, \, i_{\text{Space}} \in \mathcal{J}_{\text{Space}}(s) \right\}.
	$$
	Let us emphasize the independence of the polynomial degree in time from the degrees of the spatial polynomials. It is motivated by the fact that the complexity of the value function with respect to time might heavily differ from its complexity in space. In the linear quadratic case for example the value function is known to be quadratic in space while the behaviour with respect to time can be much more complex and is unknown, c.f.~\eqref{eq: KalVF}. In such a setting it is reasonable to allow for a higher degree in the time polynomials while choosing a small hyperbolic cross index for the spatial polynomials. Note that this approach is in line with choosing the number of spatial and time samples individually, c.f.~\Cref{subsec: data}.
	
	%%%%%%%%%%%%%%%%%%%%%%%%%%%%%%%%%%%%%%%%%%%%%%%%%%%%%%%%
	\subsection{Obtaining polynomial approximation by linear regression}\label{subsec: LeSq}~\\
	%%%%%%%%%%%%%%%%%%%%%%%%%%%%%%%%%%%%%%%%%%%%%%%%%%%%%%%%
	
	We finally recover the polynomial approximation of the value function by fitting its truncated polynomial representation to the generated data of the form \eqref{eq: sampl data}. In order to allow a comparison of the performance improvements achieved by including gradients and Hessians in the fitting process we introduce three weights $\beta_0$, $\beta_1$, $\beta_2 \geq 0$ corresponding to the values of $\mathcal{V}$, its gradient, and its Hessian respectively. 
	
	For the construction of the vector of sampled data we define
	\begin{equation*}
		V_\mathcal{V} = \frac{1}{\sqrt{N}} \left( \mathcal{V}(t_i,\xi_i) \right)_{i = 1}^N \in \mathbb{R}^N,
	\end{equation*}
	further for $j = 1,...,n$ we set
	\begin{equation*}
		V_{\nabla,j} = \frac{1}{\sqrt{N}} \left( \frac{\partial}{ \partial \xi_j} \mathcal{V}(t_i,\xi_i) \right)_{i = 1}^N \in \mathbb{R}^N,
	\end{equation*}
	and finally for $k,h = 1,...,n$ satisfying $k \leq h$ we define
	\begin{equation*}
		V_{\nabla^2,k,h} = \frac{1}{\sqrt{N}} \left( \frac{\partial^2}{ \partial \xi_k \xi_h} \mathcal{V}(t_i,\xi_i) \right)_{i = 1}^N \in \mathbb{R}^N.
	\end{equation*}
	The corresponding matrices are constructed as follows. For the function values we set
	\begin{equation*}
		A_\mathcal{V} = \frac{1}{\sqrt{N}} \left( \Phi_j(t_i,\xi_i) \right)_{i = 1,j \in \mathcal{J}}^{i = N} \in \mathbb{R}^{N \times \vert \mathcal{J} \vert},
	\end{equation*}
	further for $k=1,...,n$ we set
	\begin{equation*}
		A_{\nabla,k} = \frac{1}{\sqrt{N}} \left( \frac{\partial}{\partial \xi_k} \Phi_j(t_i,\xi_i) \right)_{i = 1,j \in \mathcal{J}}^{i = N} \in \mathbb{R}^{N \times \vert \mathcal{J} \vert},
	\end{equation*}
	and finally for $k,h = 1,...,n$ with $k \leq h$ we set
	\begin{equation*}
		A_{\nabla^2,k,h} = \frac{1}{\sqrt{N}} \left( \frac{\partial^2}{\partial \xi_k \xi_h} \Phi_j(t_i,\xi_i) \right)_{i = 1,j \in \mathcal{J}}^{i = N} \in \mathbb{R}^{N \times \vert \mathcal{J} \vert}.
	\end{equation*}
	Note that for the Hessian samples and polynomial evaluations we enforce $k \leq h$ to exploit the symmetry of the Hessian matrix. This reduces the size of the least squares problem by $N \frac{ \, n (n-1)}{2}$. We finally set
	\begin{equation*}
		\mathbf{V} = 
		\begin{pmatrix}
			\beta_0 \,  V_\mathcal{V} \\
			\beta_1 \, V_{\nabla,1} \\
			\vdots \\
			\beta_1 \, V_{\nabla,n}\\
			\beta_2 \, V_{\nabla^2,1,1}\\
			\vdots \\
			\beta_2 \, V_{\nabla^2,n,n}
		\end{pmatrix}
		~~~~~\text{ and }~~~~~
		\mathbf{A} = 
		\begin{pmatrix}
			\beta_0 \, A_\mathcal{V}\\
			\beta_1 \, A_{\nabla,1}\\
			\vdots \\
			\beta_1 \, A_{\nabla_n} \\
			\beta_2 \, A_{\nabla^2,1,1}\\
			\vdots \\
			\beta_2 \, A_{\nabla^2,n,n}
		\end{pmatrix}
	\end{equation*}
	and the coordinates $(\theta_j)_{j \in \mathcal{J}}$ of the polynomial approximation of $\mathcal{V}$ are given as the solution of the linear least squares problem
	\begin{equation}\label{eq: LeastS}
		\min_{\theta \in \mathbb{R}^{\vert \mathcal{J} \vert}} \Vert \mathbf{A} \theta - \mathbf{V} \Vert_2^2. 
	\end{equation}
	In our implementation we solve the least squares problem using the \matlab\, backslash routine, i.e., $\theta = \mathbf{A} \backslash \mathbf{V} $.
	
	\begin{remark}
		The results in our numerical experiments exhibit some minor oscillations. In particular the Hessian of the polynomial approximations are prone to this issue which might lead to noticeable deviations in the inverse. Such issues can be tackled by introducing a regularizing term in \eqref{eq: LeastS}. Note, however, that even the implementation of a simple $L^2$-penalty term is non-trivial in our setup. This is due to the augmentation of the sampling by including derivatives. Furthermore, this adjustment comes with a substantial increase in computational cost. Due to the smoothing effect of integrating the observer equation \eqref{eq: observerEq} the moderate oscillations on the inverse Hessian did not pose an immediate problem to our approach. We therefore decided against the application of a regularizer in the least squares problem.
		
		Including an appropriate regularizer can furthermore ensure sparsity of the solution. We refer to \cite{Azmetal21} where a \textit{weighted LASSO} was deployed.
	\end{remark}
	
	%%%%%%%%%%%%%%%%%%%%%%%%%%%%%%%%%%%%%%%%%%%%%%%%%%%%%%%%
	\section{Realization of the observer trajectory}\label{sec: realObs}~\\
	%%%%%%%%%%%%%%%%%%%%%%%%%%%%%%%%%%%%%%%%%%%%%%%%%%%%%%%%
	
	After obtaining a polynomial approximation of the value function $\mathcal{V}_\mathrm{p}$ we are in the position to realize the Mortensen observer trajectory $\widehat{x}_\mathrm{M}$ numerically. As announced this is done using two different approaches namely by solving the observer equation \eqref{eq: observerEq} and by minimizing the value function according to \eqref{eq: observerDef}.
	
	\vspace{5mm}
	
	\textbf{Minimizing the value function.}
	We construct a a discrete approximation of the Mortensen observer trajectory on a time grid of $[0,T]$ using $N_{\mathrm{min}} + 1$ equidistant grid points, where $t_0 = 0$ and $t_{N_\mathrm{min}} = T$. For any $k \in \{ 0,...,N_\mathrm{min} \}$ we construct the observer trajectory by setting
	\begin{equation*}
		\widehat{x}_{\mathrm{min},k} = \widehat{x}_\mathrm{min}(t_k) = \argmin_{\xi \in \mathbb{R}^n} \mathcal{V}_\mathrm{p}(t_k,\xi).
	\end{equation*}
	The minimization is realized using a gradient scheme with an \textit{Armijo} stepsize control. For the construction of $\widehat{x}_{\mathrm{min},0}$ we initialize the minimization by $x_0$, which is the exact minimizer of $\mathcal{V}(0,\cdot)$. For the construction of $\widehat{x}_{\mathrm{min},k}$ with $k > 0$ the initialization for the minimization is set to $\widehat{x}_{\mathrm{min},k-1}$. Assuming that $\widehat{x}_\mathrm{M}$ is continuous and that $N_\mathrm{min}$ was chosen large enough, this yields a sufficiently good initial estimate justifying the choice made for the sampling points in \Cref{subsec: data}.
	
	\vspace{5mm}
	
	\textbf{Solving the observer equation.}
	Another approximation of the observer trajectory $\widehat{x}_\mathrm{eq}$ is computed by solving the observer equation where the inverse Hessian of the value function is approximated by $\nabla_{\xi\xi}^2 \mathcal{V}_\mathrm{p}(\cdot,\cdot)^{-1}$. Here a grid of $N_\mathrm{eq} + 1$ points with $t_0 = 0$ and $t_{N_\mathrm{eq}} = T$ is used. The equation is solved using a BDF4 scheme, where non-linearities are handled using a Newton scheme.
	
	%%%%%%%%%%%%%%%%%%%%%%%%%%%%%%%%%%%%%%%%%%%%%%%%%%%%%%%%%%%%%%
	\section{Numerical Tests}\label{sec: num}
	%%%%%%%%%%%%%%%%%%%%%%%%%%%%%%%%%%%%%%%%%%%%%%%%%%%%%%%%%%%%%%
	
	In this section we apply the proposed methodology to four different models. First we consider a linear model allowing a comparison with the Kalman filter. We further illustrate the benefits of including the sampled Hessian in the regression problem and compare the realizations via a minimization of the value function according to \eqref{eq: observerDef} and via solving the observer equation \eqref{eq: observerEq}. In the second and third example  we consider two non-linear low-dimensional oscillators and set the focus on comparing the Mortensen observer with the extended Kalman filter. These examples further illustrate the challenges stemming from the fact that the systems need to be considered backwards in time. As a final example we present an agent based model of higher state-space dimension.
	
	In the linear case the Kalman filter and the Mortensen observer are equivalent on a theoretical level. Therefore we have access to an accurate approximation of the Mortensen observer that our numerical results can be compared to. Since this is not the case for the non-linear problems, we instead solve the observer equation \eqref{eq: observerEq} using a BDF4 scheme where we compute the inverse of the Hessian $\nabla_{\xi\xi}^2 \mathcal{V}(t,x)^{-1}$ by solving the corresponding open loop problem before solving the DRE as described in \Cref{subsec: data}. The resulting trajectory will be denoted by $\widehat{x}_{\mathrm{M}}$ and we will consider it to be the true Mortensen observer trajectory.
	
	%%%%%%%%%%%%%%%%%%%%%%%%%%%%%%%%%%%%%%%%%%%%%%%%%%%%%%%%%%%%%%
	\subsection{Practical aspects}\label{subs: practAsp}~\\
	%%%%%%%%%%%%%%%%%%%%%%%%%%%%%%%%%%%%%%%%%%%%%%%%%%%%%%%%%%%%%%
	
	We first turn our attention to the evaluation of the value function and its derivatives in the sampling points. Since these calculations are entirely independent of each other they can easily be parallelized. In our implementation this is done using \texttt{parfor} from the \matlab\, Parallel Computing Toolbox.
	
	The gradient descent scheme is implemented with a relative tolerance of $10^{-6}$. For the non-linear examples we further implement an absolute tolerance of $10^{-3}$ and terminate the scheme once one of the tolerances is reached. The state and adjoint equations are solved by an application of a BDF4 scheme using 1001 time discretization points. The non-linearities are treated by a Newton scheme with absolute tolerance $10^{-12}$. The Riccati equations are solved using the same BDF4 scheme where the implicit time steps are realized by the \matlab\, routine \texttt{icare}. Only for the four initial time steps a Newton scheme is employed, where the absolute tolerance is set to $10^{-10}$ for the first three examples and to $10^{-8}$ for the fourth example.
	
	For the integration of the observer equation by means of the approximated value function we set $N_\mathrm{eq} = 10^3$ and deploy a BDF4 scheme in which non-linearities are treated using a Newton scheme with absolute tolerance $10^{-8}$. The minimization of the approximated value function is performed via a gradient descent scheme with relative tolerance $10^{-6}$ and absolute tolerance $10^{-3}$ and with $N_\mathrm{min} = 10^3$. 
	
	The trajectory $\widehat{x}_\mathrm{K}$ resulting from the (extended) Kalman filter and the solution $\Sigma$ of the corresponding Riccati equation are determined using the \matlab\, routine \texttt{ode15s} with $1001$ equidistant discretization points and a relative tolerance of $10^{-8}$. When solving for the true Mortensen trajectory via BDF4 non-linearities are treated by a Newton scheme with absolute tolerance $10^{-6}$.
	
	The schemes were implemented in \matlab\, R2020b and the computations were run on a Lenovo ThinkPad T14s AMD Ryzen 7 PRO 4750U(16)@1.700GHz with 32GB DDR4 3200 MHz memory.
	
	The \matlab\, code used to obtain the numerical results is available in \cite{SchB23}.
	
	%%%%%%%%%%%%%%%%%%%%%%%%%%%%%%%%%%%%%%%%%%%%%%%%%%%%%%%%%%%%%%
	\subsection{Test 1: Harmonic oscillator}\label{subs: test 1}~\\
	%%%%%%%%%%%%%%%%%%%%%%%%%%%%%%%%%%%%%%%%%%%%%%%%%%%%%%%%%%%%%%
	
	For a linear example we consider an undamped harmonic oscillator. The first order form of the model reads
	\begin{equation*}
		\begin{aligned}
			\dot{x}(t) &= A x(t) + F v(t),~
			x(0) = x_0,\\
			y(t) &= C x(t) + \mu(t),
		\end{aligned}
	\end{equation*}
	where
	$$A = \begin{pmatrix} 0 & 1\\ -1 & 0 \end{pmatrix}, F = \begin{pmatrix} 0 \\ 1 \end{pmatrix}, C = \begin{pmatrix} 1 & 0 \end{pmatrix}, x_0 = \begin{pmatrix}  1 \\ 1 \end{pmatrix}.$$
	In this formulation the first component of the state represents the position of the oscillator while the second component corresponds to its velocity. In our computations the disturbance of the dynamics is restricted to the velocity while the observation measures only the position. For the artificial construction of the measured data $y$ we set the error in the dynamics to  $v(t) = \tfrac{1}{2} \cos(\tfrac{6}{5}t)$ and the observation error to $\mu(t) = \tfrac{1}{2} \sin(\tfrac{t}{2})$.
	
	In order to reconstruct the state we apply our previously described methodology to approximate the Mortensen observer. In these calculations we consider the time horizon $[0,20]$. For the polynomial approximation of the value function we set the hyperbolic cross index to $s=5$. Note that for this linear example we know a priori that the value function is quadratic in the spatial variable. For an accurate and cost effective approximation of $\mathcal{V}$ one would restrict the spatial polynomials to a maximum degree of two. For general non-linear examples, however, such a priori knowledge is not available. We therefore decided to include spatial polynomials of higher degrees in the basis. For each computation the maximum degree for the time polynomials is set to be equal to the number of time samples used. The domains $Q^k$ introduced in \eqref{eq: samp rec} for the spatial sampling are defined via the side lengths $r_{k,1} = r_{k,2} = \max \left\{  0.1, 0.1 \, \Vert \widehat{x}_\mathrm{K}(t_k) \Vert \right\} $.
	
	In \Cref{table: Hess err} we report results for different choices of sample numbers and sampling strategies (represented by the weights $(\beta_0,\beta_1,\beta_2)$): (i) The weights $(1,0,0)$ corresponds to the classical scheme of only considering the value function values. (ii) With the weights $(10^{-3},1,0)$ we pay attention mostly to the gradients. (iii) The combination $(10^{-3},0,1)$ implies a focus on the Hessian. (iv) Finally the weights $(1,1,\tfrac{1}{2})$ are derived from the coefficients of the general Taylor polynomial and imply a consideration of all available information. Note that the different weights result in different sizes of the least squares problem. Neglecting the computational effort of evaluating the samples the burden of approximating the value function lies in the least squares problem. We therefore decided to compare the estimated computational cost resulting from the different strategies by comparing the number of rows of the least squares matrix and display it in the fourth column of \Cref{table: Hess err}. In the last three columns we present the relative errors associated with the approximations of the Mortensen observer. The sixth column displays the relative error of the trajectory $\widehat{x}_\mathrm{min}$ obtained via minimization of the approximated value function $\mathcal{V}_\mathrm{p}$, i.e.,
	\begin{equation*}
		e_\mathrm{min} = \frac{\Vert \widehat{x}_\mathrm{min} - \widehat{x}_\mathrm{K} \Vert_{L^2(0,T;\mathbb{R}^n)}}{\Vert \widehat{x}_\mathrm{K} \Vert_{L^2(0,T;\mathbb{R}^n)}}.
	\end{equation*}
	The seventh column shows the relative error of the observer trajectory $\widehat{x}_\mathrm{eq}$ obtained by solving the observer equation using the Hessian of the approximation of the value function $\nabla_{\xi\xi}^2\mathcal{V}_\mathrm{p}$, i.e.,
	\begin{equation*}
		e_\mathrm{eq} = \frac{\Vert \widehat{x}_\mathrm{eq} - \widehat{x}_\mathrm{K} \Vert_{L^2(0,T;\mathbb{R}^n)}}{\Vert \widehat{x}_\mathrm{K} \Vert_{L^2(0,T;\mathbb{R}^n)}}.
	\end{equation*}
	We further present the relative error of the inverted Hessian
	\begin{equation*}
		e_\mathrm{gain} =
		\frac{\Vert \nabla_{\xi\xi}^2 \mathcal{V}_\mathrm{p}(\cdot,\widehat{x}_\mathrm{M}(\cdot))^{-1} - \Sigma(\cdot) \Vert_{L^2(0,T;\mathbb{R}^{n,n})}}
		{\Vert \Sigma \Vert_{L^2(0,T;\mathbb{R}^{n,n})}}.
	\end{equation*}
	All entries of the table where $e_\mathrm{min}$ or $e_\mathrm{eq}$ are marked as \texttt{failed} represent parameter combinations that resulted in a polynomial approximation $\mathcal{V}_\mathrm{p}$ for which the minimization or respectively the integration of the observer equation did not converge. 
	
	Our experiments suggest that including the Hessian samples in the least squares problems leads to higher accuracy both in the approximated inverted Hessian and in the resulting observer trajectories. In row 11 the least squares matrix considering the Hessian information has 600 rows and yields a polynomial approximation of the value function based on which we approximate the Mortensen observer via integration of \eqref{eq: observerEq} with a relative error of order $10^{-6}$. In our experiments we did not reach this level of accuracy without including the sampled Hessians in the least squares problem. Row 2 of \Cref{table: Hess err} shows that even an increase of the number of rows in the LS-problem to $1800$ results in an approximation of the observer equation with a relative error of order $10^{-3}$. Comparing, e.g., row 11 and row 12 we further learn that sampling only the values and the Hessians seems to be more effective than additionally including the gradients. Not only do the gradients increase the size of the LS-problem, in most cases the resulting approximation of the value function also yields less accurate approximations of the Mortensen observer. Only in row 19 and 20 the two combinations of weights lead to observer trajectories with the same level of accuracy when integrating the observer equation.
	
	We further observe that in all parameter constellations the integration of the observer equation yields results preferable to the ones obtained by minimization of the value function. In row 6 we even find an example for which the observer equation leads to an approximation with a relative error of order $10^{-2}$ while the minimization of the value function fails. 
	
\begin{figure}
		\centering
		\includegraphics[scale = 0.8]{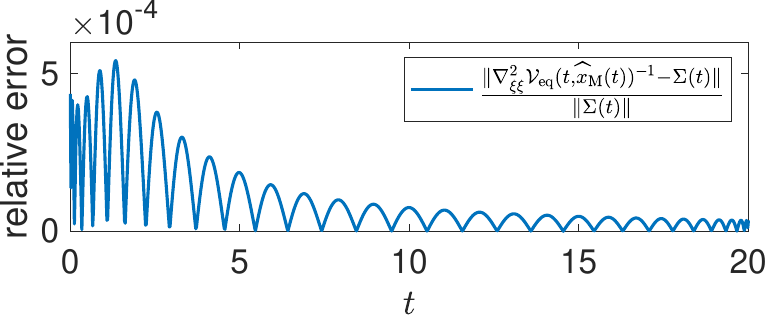}
		\caption{Relative error of the inverted Hessian measured in the Frobenius norm plotted along time. The computation is based on a polynomial approximation of $\mathcal{V}$ obtained with the parameters $N_{\mathrm{Time}} = 30$, $N_{\mathrm{Space}} = 5$, $(\beta_0,\beta_1,\beta_2) = (10^{-3},0,1)$, $d_{\mathrm{Time}} = 30$, and $s = 5$.}
		\label{fig: feedback error}
	\end{figure}

	To conclude the discussion of the linear example we present \Cref{fig: feedback error} illustrating the development of the relative error of the inverted Hessian along time corresponding to the parameters set in row 11 of \Cref{table: Hess err}.
	
	\begin{table}
		\begin{center}
			\begin{NiceTabular}{c | | c c c c | c c c}
				\hline
				& & & & & & relative errors  \\   
				\hline
				& $N_{\mathrm{Time}}$ & $N_{\mathrm{Space}}$ & $(\beta_0,\beta_1,\beta_2)$ & rows & $e_{\mathrm{gain}}$ & $e_{\mathrm{min}}$ & $e_{\mathrm{eq}}$ \\
				\hline \hline 
				1 & $30$ & $20$ & $(1,0,0)$ & $600$ & $1.7 \times 10^{-4}$  & $2.5 \times 10^{-4}$ & $\mathbf{1.0 \times 10^{-5}}$  \\
				2 & $30$ & $20$ & $(10^{-3},1,0)$ & $1800$ & $4.7 \times 10^{-2}$ & $2.7 \times 10^{-2}$ & $\mathbf{2.3 \times 10^{-3}}$ \\
				3 & $30$ & $20$ & $\mathbf{(10^{-3},0,1)}$ & $\mathbf{2400}$ & $1.5 \times 10^{-4}$ & $2.4 \times 10^{-4}$ & $\mathbf{7.0 \times 10^{-6}}$ \\
				4 & $30$ & $20$ & $(1,1,0.5)$ & $3600$ & $7.4 \times 10^{-4}$ & $5.2 \times 10^{-4}$ & $\mathbf{5.3 \times 10^{-5}}$ \\
				\hline
				5 & $30$ & $10$ & $(1,0,0)$ & $300$ & $18.2$ & failed & failed \\
				6 & $30$ & $10$ & $(10^{-3},1,0)$ & $900$ & $1.5 \times 10^{-1}$ & failed & $\mathbf{1.4 \times 10^{-2}}$ \\
				7 & $30$ & $10$ & $\mathbf{(10^{-3},0,1)}$ & $\mathbf{1200}$ & $1.5 \times 10^{-4}$ & $2.4 \times 10^{-4}$ & $\mathbf{7.0 \times 10^{-6}}$ \\
				8 & $30$ & $10$ & $(1,1,0.5)$ & $1800$ & $8.7 \times 10^{-4}$ & $7.1 \times 10^{-4}$ & $\mathbf{6.6 \times 10^{-5}}$ \\
				\hline
				9 & $30$ & $5$ & $(1,0,0)$ & $150$ & $86.4$ & failed & failed \\
				10 & $30$ & $5$ & $(10^{-3},1,0)$ & $450$ & $20.4$ & failed & failed \\
				11 & $30$ & $5$ & $\mathbf{(10^{-3},0,1)}$ & $\mathbf{600}$ & $1.5 \times 10^{-4}$ & $2.4 \times 10^{-4}$ & $\mathbf{7.0 \times 10^{-6}}$ \\
				12 & $30$ & $5$ & $(1,1,0.5)$ & $900$ & $3.5 \times 10^{-3}$ & $2.2 \times 10^{-3}$ & $\mathbf{9.7 \times 10^{-5}}$ \\
				\hline
				13 & $20$ & $5$ & $(1,0,0)$ & $100$ & $8.6$ & failed & failed \\
				14 & $20$ & $5$ & $(10^{-3},1,0)$ & $300$ & $24.0$ & failed & failed \\
				15 & $20$ & $5$ & $\mathbf{(10^{-3},0,1)}$ & $\mathbf{400}$ & $2.4 \times 10^{-3}$ & $2.4 \times 10^{-3}$ & $\mathbf{2.3 \times 10^{-4}}$ \\
				16 & $20$ & $5$ & $(1,1,0.5)$ & $600$ & $7.9 \times 10^{-3}$ & $5.6 \times 10^{-3}$ & $\mathbf{3.2 \times 10^{-4}}$ \\
				\hline
				17 & $10$ & $5$ & $(1,0,0)$ & $50$ & $30.3$ & failed & failed \\
				18 & $10$ & $5$ & $(10^{-3},1,0)$ & $150$ & $6.6$ & failed & failed \\
				19 & $10$ & $5$ & $\mathbf{(10^{-3},0,1)}$ & $\mathbf{200}$ & $3.8 \times 10^{-2}$ & $6.4 \times 10^{-1}$ & $\mathbf{3.5 \times 10^{-3}}$ \\
				20 & $10$ & $5$ & $(1,1,0.5)$ & $300$ & $3.8 \times 10^{-2}$ & $6.4 \times 10^{-1}$ & $\mathbf{3.5 \times 10^{-3}}$ \\
				\hline
			\end{NiceTabular}
			\caption{Comparison of the relative errors of the inverted Hessian and the two observer trajectories obtained by minimization of $\mathcal{V}$ and by solving the observer equation depending on the number of samples and weights used in the least squares problem.}
			\label{table: Hess err}
		\end{center}
	\end{table}
	
	%%%%%%%%%%%%%%%%%%%%%%%%%%%%%%%%%%%%%%%%%%%%%%%%%%%%%%%%%%%%%%
	\subsection{Test 2: Van der Pol oscillator}\label{subs: test 2}~\\
	%%%%%%%%%%%%%%%%%%%%%%%%%%%%%%%%%%%%%%%%%%%%%%%%%%%%%%%%%%%%%%
	
	We now turn our attention to the first non-linear example. Again we approximate the Mortensen observer both via the minimization of $\mathcal{V}_\mathrm{p}$ and via integration of the observer equation. Both approaches yield satisfying results but again the latter requires less data to do so. We further compare our findings to the results obtained using the extended Kalman filter. 
	
	Specifically we consider the van der Pol oscillator. It is modelled by 
	\begin{equation*}
		\begin{aligned}
			\frac{\mathrm{d}}{\mathrm{d}t} \begin{pmatrix} x_1(t) \\ x_2(t) \end{pmatrix}  
			&= A \begin{pmatrix} x_1(t) \\ x_2(t) \end{pmatrix} 
			- \begin{pmatrix} 0 \\ x_1(t)^2 ~ x_2(t) \end{pmatrix}
			+ F v(t),~
			\begin{pmatrix} x_1(0) \\ x_2(0) \end{pmatrix} = x_0,\\
			y(t) &= C \begin{pmatrix} x_1(t) \\ x_2(t) \end{pmatrix} + \mu(t),
		\end{aligned}
	\end{equation*}
	where
	\begin{equation*}
		A = \begin{pmatrix} 0 & 1\\ -1 & 1 \end{pmatrix},~ F = \begin{pmatrix} 0 \\ 1 \end{pmatrix},~ C = \begin{pmatrix} 1 & 0 \end{pmatrix},~ x_0 = \begin{pmatrix} 0.1 \\ 0.1 \end{pmatrix}.
	\end{equation*}
	The two state variables $x_1$ and $x_2$ correspond to the position and velocity of the system, respectively. Note that by introducing a third variable $x_3 = x_1^2$ this system can be transformed into one of state dimension $n = 3$ with a quadratic right hand side. Therefore theoretical results presented in \cite{BreS23} apply.
	
	The behaviour of this system is best characterized in terms of the phase space. All initial values are attracted to a limit cycle. Once the limit cycle is reached the state of the system will stay near that cycle indefinitely. Hence this system exhibits stable behaviour. However, when considered backwards in time the system turns unstable. Especially when starting in a point outside the limit cycle and observing the evolution of the system backwards in time the norm of the state will increase rapidly. These dynamics pose a crucial issue in our approach of realizing the Mortensen observer, specifically with regard to solving the open loop control problem for a given sample point $(t^*,\xi^*)$. In particular for large $t^*$ and $\xi^*$ close to the limit cycle the gradient descent scheme requires an accurate initial guess for the optimal control. In our experiments we tackled this issue in an iterative manner. For a fixed integer $N_\mathrm{it}$ we first solve the optimal control problem for the tuple $(\frac{t^*}{N_\mathrm{it}},\xi^*)$ using the zero control for initialization. The resulting optimal control is used as an initialization when solving the open loop problem for $(\frac{2 \,t^*}{N_\mathrm{it}},\xi^*)$. Therefore we have to solve $N_\mathrm{it}$ control problems in order to obtain the desired evaluation of the sampling point.
	
	For the construction of the measured data $y$ we set the disturbance in the dynamics to $v(t) = \tfrac{1}{2} \cos(\tfrac{6}{5} \, t)$ and the disturbance in the observation is chosen as $\mu(t) = \tfrac{3}{10} \sin(2 \pi \, t)$. This example is considered over the time horizon $[0,7]$ and the domains for spatial sampling are set via the side lengths $r_{k,1} = r_{k,2} = \max \left\{  0.1, 0.1 \, \Vert \widehat{x}_\mathrm{K}(t_k) \Vert \right\} $.
	
\begin{figure}
		\centering
		\begin{subfigure}{\textwidth}
			\includegraphics[scale = 0.295]{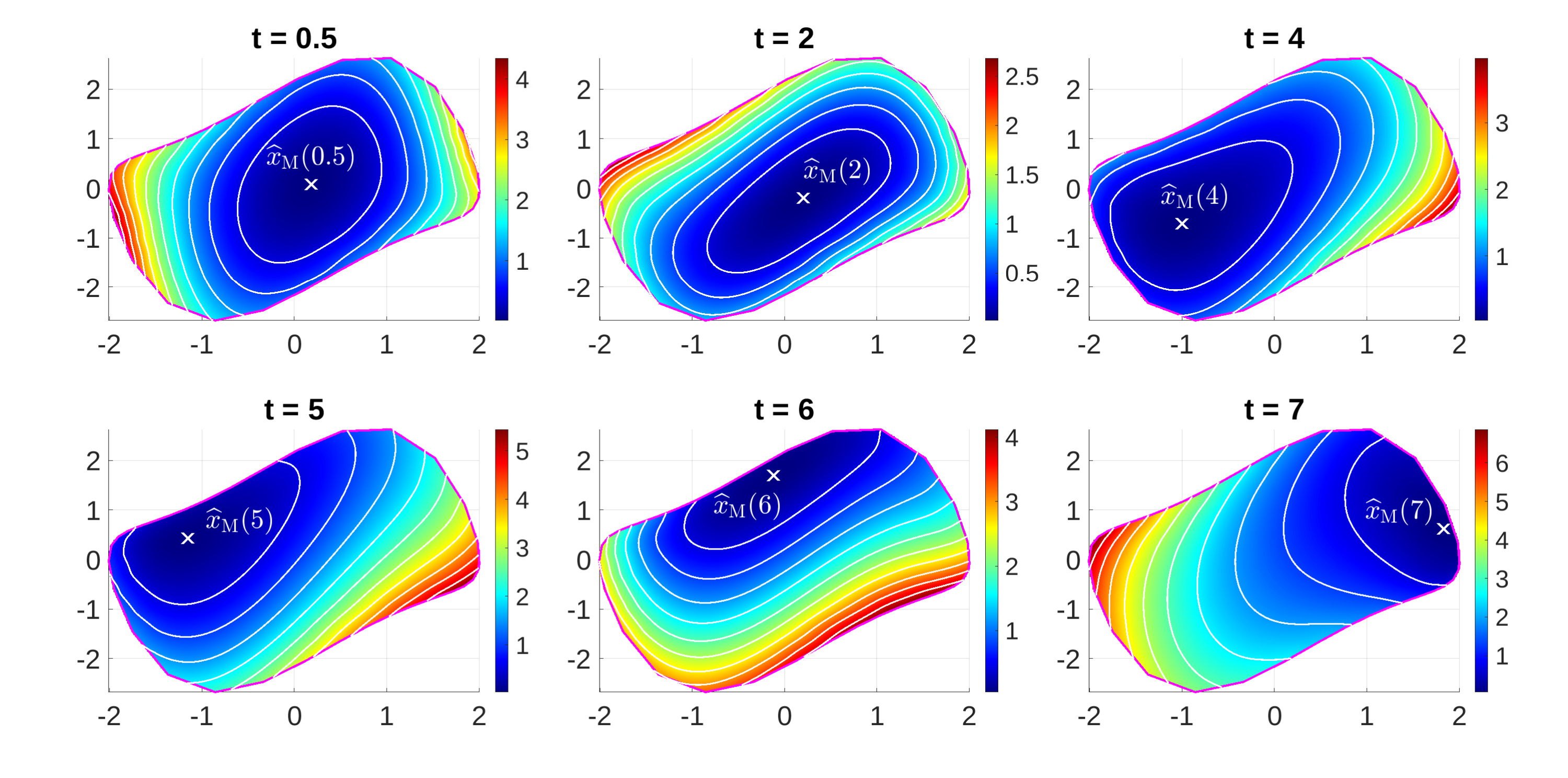}
			\caption{Vale function for fixed times $t$} 
			\label{fig: VdP_VF}
		\end{subfigure}
		\begin{subfigure}{0.45\textwidth}
			\includegraphics[scale = 0.5]{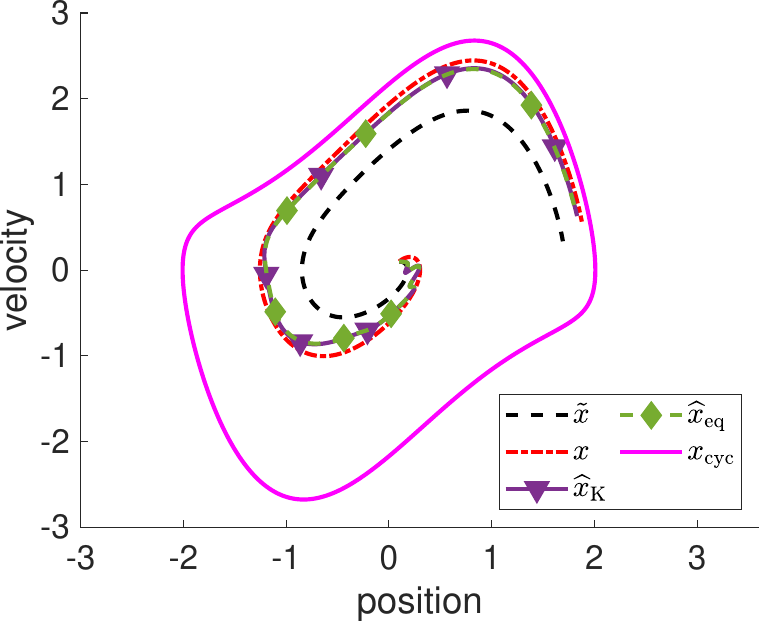}
			\caption{Approximation by solving observer equation with pol. appr. using 750 samples}
			\label{fig: VdP_eq}
		\end{subfigure}
		\hfill
		\begin{subfigure}{0.45\textwidth}
			\includegraphics[scale = 0.5]{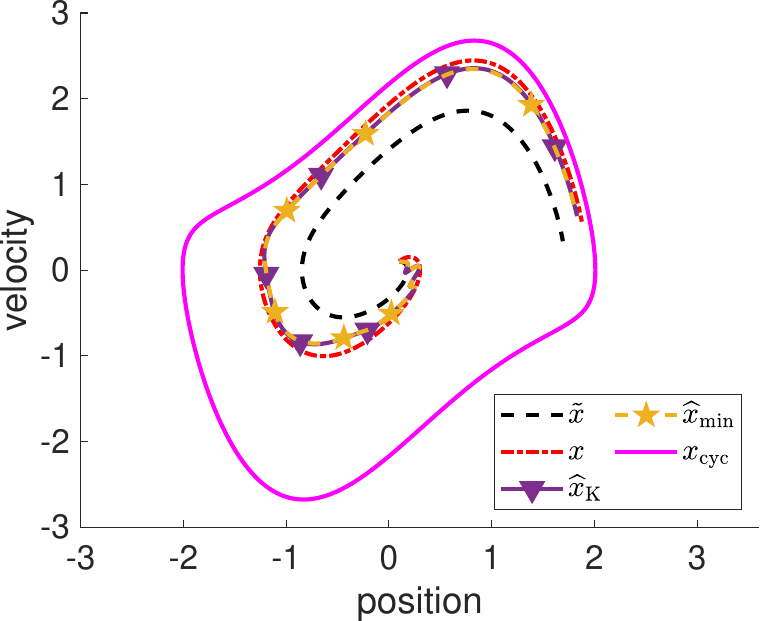}
			\caption{Approximation by minimizing pol. appr. of the value function using 3000 samples}
			\label{fig: VdP_min}
		\end{subfigure}
		\caption{Van der Pol oscillator}
		\label{fig: VdP}
	\end{figure}

	The results are illustrated in \Cref{fig: VdP}. In \Cref{fig: VdP_VF} we show a plot of the value function $\mathcal{V}$ evaluated for fixed times $t$ for $\xi$ inside the limit cycle.
	
	In order to integrate the observer equation the value function was approximated using $30$ and $25$ sampling points in time and space, respectively. In the LS-problem only the values and the Hessians where considered with respective weights $10^{-3}$ and $1$ resulting in a LS-matrix with $3000$ rows. The polynomial basis is set via $d_{Time} = 9$ and $s = 9$. The computation took roughly 20 minutes and the results are presented in \Cref{fig: VdP_eq}. The relative error of the obtained trajectory is given by
	\begin{equation*}
		\frac{\Vert \widehat{x}_\mathrm{eq} - \widehat{x}_\mathrm{M} \Vert_{L^2(0,7;\mathbb{R}^2)}}{\Vert \widehat{x}_\mathrm{M} \Vert_{L^2(0,7;\mathbb{R}^2)}}
		= 1.8 \times 10^{-3}.
	\end{equation*}
	
	For the realization via the minimization of $\mathcal{V}_\mathrm{p}$ we used $60$ and $50$ time an space samples, respectively. In the LS-problem only values and gradients were included using weights $10^{-3}$ and $1$, respectively, hence the LS-matrix has $9000$ rows. The polynomial basis is given by $d_\mathrm{Time} = 17$ and $s = 10$. Here the computation took about 90 minutes and the results are presented in \Cref{fig: VdP_min}. They exhibit a relative error of
	\begin{equation*}
		\frac{\Vert \widehat{x}_\mathrm{min} - \widehat{x}_\mathrm{M} \Vert_{L^2(0,7;\mathbb{R}^2)}}{\Vert \widehat{x}_\mathrm{M} \Vert_{L^2(0,7;\mathbb{R}^2)}}
		= 3.3 \times 10^{-3}.
	\end{equation*}
	
	From \Cref{fig: VdP} we observe that for this particular system the Mortensen observer and the extended Kalman filter lead to very similar trajectories for the reconstruction of the state. The following example provides a situation where such similarities do not occur. 
	
	%%%%%%%%%%%%%%%%%%%%%%%%%%%%%%%%%%%%%%%%%%%%%%%%%%%%%%%%%%%%%%
	\subsection{Test 3: Duffing oscillator}\label{subs: test 3}~\\
	%%%%%%%%%%%%%%%%%%%%%%%%%%%%%%%%%%%%%%%%%%%%%%%%%%%%%%%%%%%%%%
	
	The purpose of the following example is to show that for more complex systems there is a substantial difference between the extended Kalman filter and the Mortensen observer. To this end we consider the Duffing equation which in its state-space form is given by
	\begin{equation}\label{eq: Duffing}
		\begin{aligned}
			\frac{\mathrm{d}}{\mathrm{d}t} \begin{pmatrix} x_1(t) \\ x_2(t) \end{pmatrix} &= 
			A
			\begin{pmatrix} x_1(t) \\ x_2(t) \end{pmatrix} 
			+ 
			\begin{pmatrix}
				0 \\ 
				- \beta \, x_1(t)^3
			\end{pmatrix} + F \, v(t),~~~
			\begin{pmatrix} x_1(0) \\ x_2(0) \end{pmatrix} = x_0 + \eta,\\
			y(t) &= C \begin{pmatrix} x_1(t) \\ x_2(t) \end{pmatrix} + \mu(t),
		\end{aligned}
	\end{equation}
	where
	\begin{equation*}
		A= \begin{pmatrix}
			0 & 1 \\
			- \lambda & - \delta
		\end{pmatrix},~
		F = \begin{pmatrix} 0 \\ 1 \end{pmatrix},~
		C = \begin{pmatrix} 1 & 0\end{pmatrix}.
	\end{equation*}
	
	A thorough discussion of equations of this type can be found in \cite{JoSm07}. Again we point out that by introducing a third variable as $x_3 = x_1^2$ this system is equivalent to one of dimension $n = 3$ with a quadratic right hand side and results from \cite{BreKu21} can be applied. However, the results presented there require the assumption that the difference of modeled and measured output $\Vert y - C \Tilde{x} \Vert_{L^2(0,T;\mathbb{R}^r)}$ is sufficiently small. Here $\Tilde{x}$ is the model trajectory, i.e., the solution of the undisturbed model equation. For systems as sensitive as the Duffing oscillator this is a rather strong assumption because even small disturbances in the dynamics and in the initial value may cause major differences in the resulting trajectory. 
	
	Just like the previous example this system is affected by the issue of backwards instability. In order to compute the evaluations in the sampling points we apply the iteration described in \Cref{subs: test 2}.
	
	We construct the measurement $y$ as follows:~In order to prompt chaotic behaviour of the system from $t = 0$ onwards we follow \cite{JoSm07,BreKu21} and set
	\begin{equation*}
		\lambda = -1,~ \beta = 1,~ \delta = 0.3,~ v(t) = \gamma \cos(\omega t),~ \gamma = 0.5,~ \omega = 1.2.
	\end{equation*}
	As the disturbance in the observation we consider $\mu(t) = 0.05 \sin(2 \pi t)$ and for the initial value and disturbance we set $x_0 = \begin{pmatrix}
		0 & 0
	\end{pmatrix}^\top$ and $\eta = \begin{pmatrix}
		-1.216 & 0.493
	\end{pmatrix}^\top$, respectively.
	
\begin{figure}
		\centering
		\begin{subfigure}{\textwidth}
			\includegraphics[scale = 0.29]{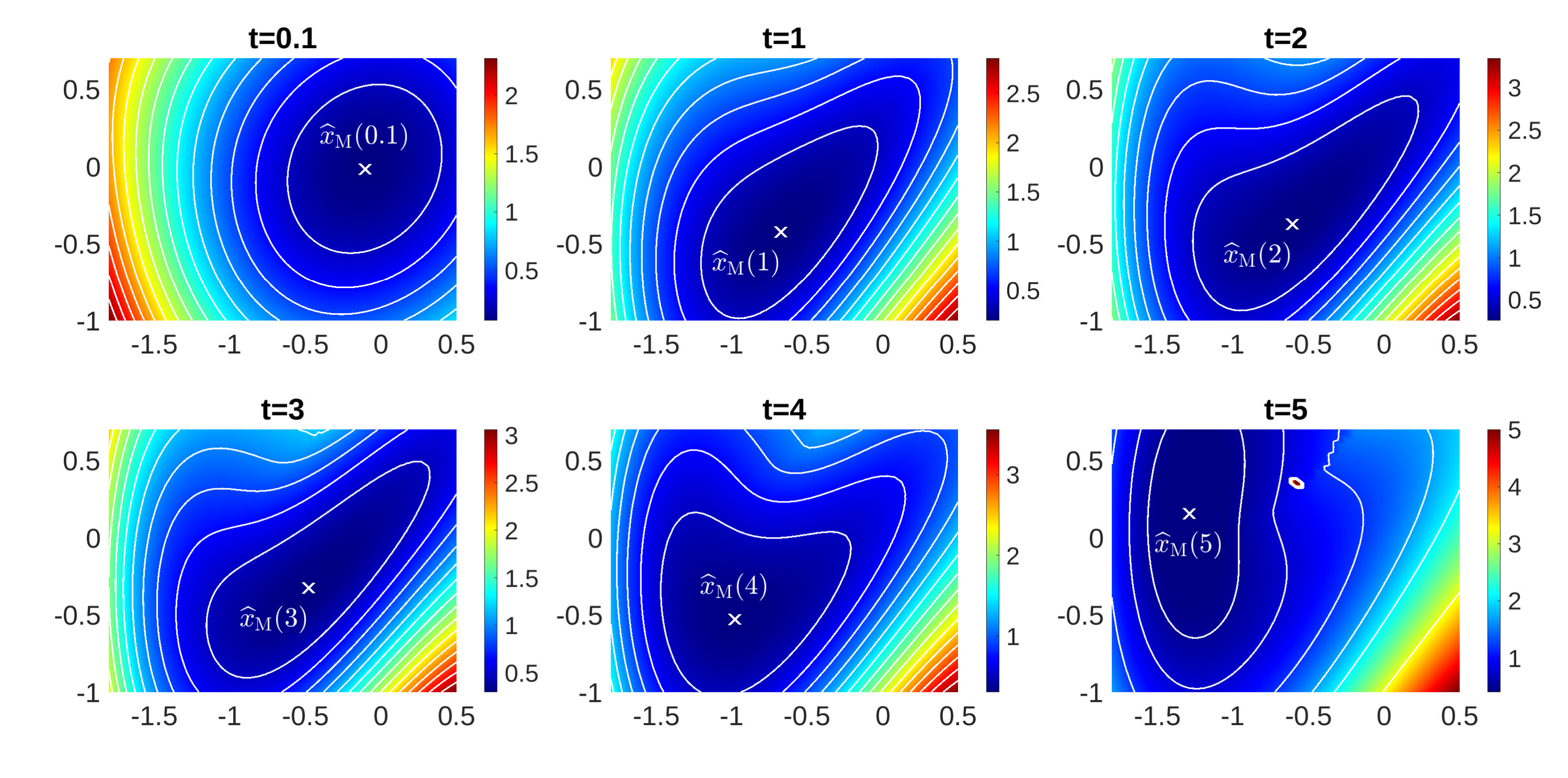}
			\caption{Value function around the minimizer}
			\label{fig: Duffing_VF}
		\end{subfigure}
		\begin{subfigure}{\textwidth}
			\includegraphics[scale = 0.5]{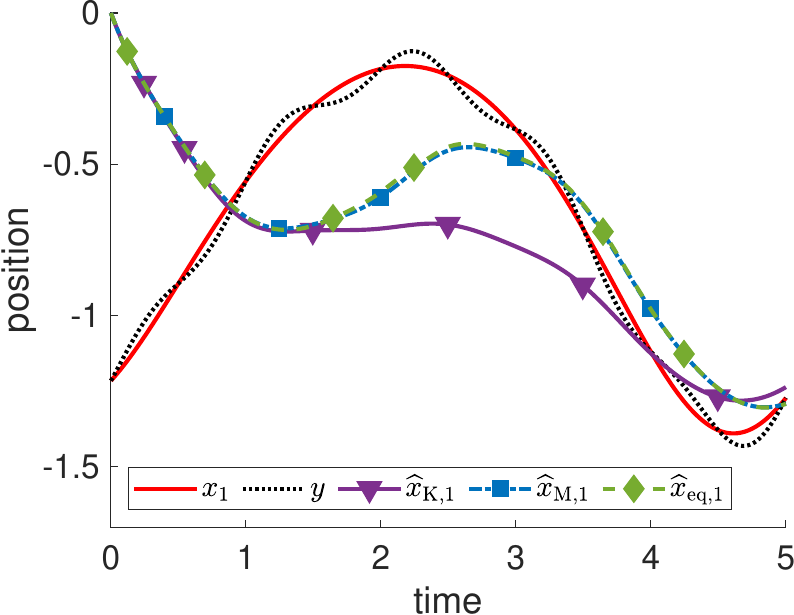}
			\hfill
			\includegraphics[scale = 0.5]{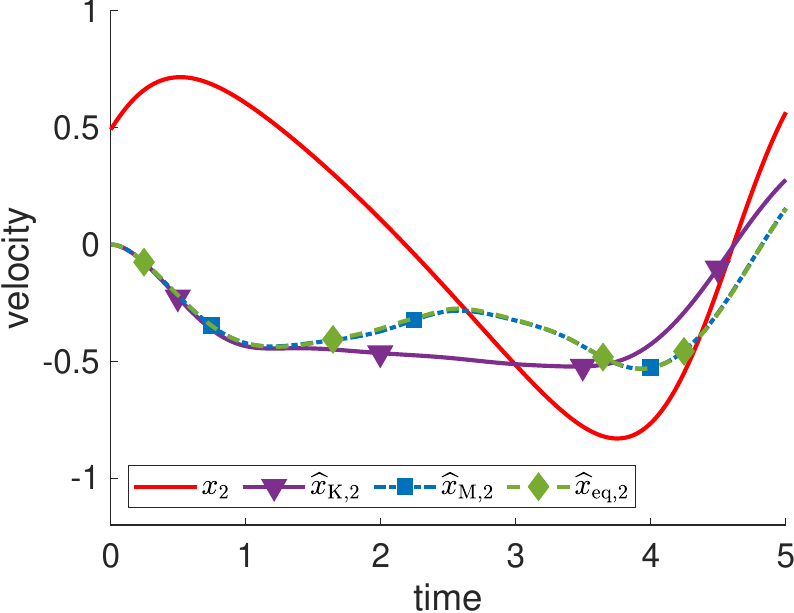}
			\caption{Approximation by solving observer equation with pol. appr. using 1050 samples }
			\label{fig: Duffing_obs}
		\end{subfigure}
		\caption{Duffing equation}
		\label{fig: Duffing}
	\end{figure}

	The results are shown in \Cref{fig: Duffing}. \Cref{fig: Duffing_VF} shows plots of the value function for fixed times $t$ evaluated around the minimizer. The frame for time $t = 5$ displays some numerical inaccuracies underlining the fact that the value function evaluation is not a trivial task. We suspect that for the final value 
	$\xi = \begin{bmatrix}
		-0.58 && 0.35
	\end{bmatrix}^\top$
	the gradient descent solving the open loop problem converged to a local instead of the global minimizer. Since this point lies outside the sampling rectangle used in our approximation scheme, this did not pose any issues while approximating the Mortensen observer.
	
	For this example we omit the minimization of $\mathcal{V}_\mathrm{p}$ and focus on the integration of the observer equation \eqref{eq: observerEq}. We consider the time horizon $[0,5]$. Here the domains for spatial sampling are defined via the side lengths $r_{k,1} = r_{k,2} = \max \left\{  0.1, 0.4 \, \Vert \widehat{x}_\mathrm{K}(t_k) \Vert \right\} $. We decided to take $35$ time and $30$ space samples. The gradients $\nabla_\xi \mathcal{V}$ will be omitted in the LS-problem while the values and Hessians are considered with the weights $10^{-2}$ and $1$, respectively. The polynomial basis is constructed using $d_\mathrm{Time} = 9$ and $s = 17$. We compare our obtained state reconstruction $\widehat{x}_\mathrm{eq}$ to the trajectory $\widehat{x}_\mathrm{K}$ obtained by means of the extended Kalman filter and with $\widehat{x}_\mathrm{M}$ constructed as was described above. The results are presented in \Cref{fig: Duffing_obs}. The reader can observe a substantial difference between the reconstructions based on the Mortensen observer and the extended Kalman filter in the case of the Duffing equation.  We also report that our approximation of the Mortensen observer has a relative error of 
	\begin{equation*}
		\frac{\Vert \widehat{x}_\mathrm{eq} - \widehat{x}_\mathrm{M} \Vert_{L^2(0,5;\mathbb{R}^2)}}{\Vert \widehat{x}_\mathrm{M} \Vert_{L^2(0,5;\mathbb{R}^2)}} = 9.4 \times 10^{-3}.
	\end{equation*}
	
\begin{figure}
		\centering
		\begin{subfigure}{0.49\textwidth}
			\includegraphics[scale = 0.57]{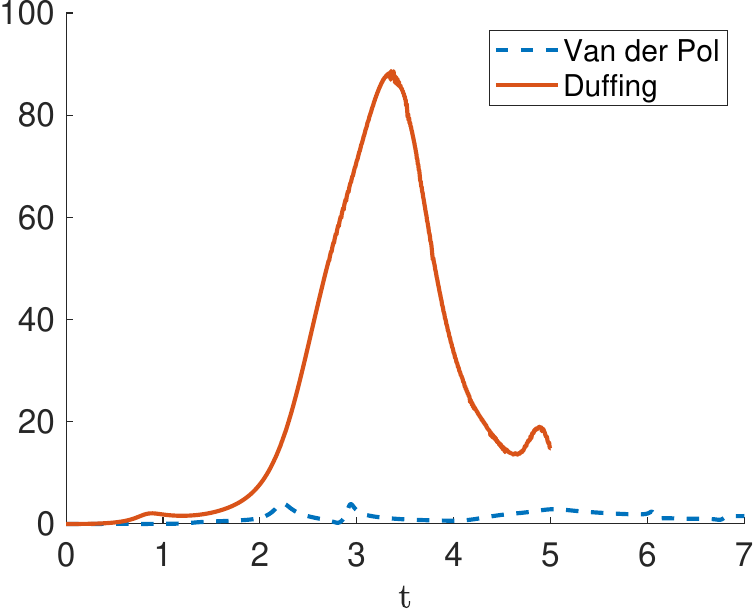}
			\caption{Relative differences of the Riccati type equations}
			\label{fig: DiffThirdDer}
		\end{subfigure}
		\hfill
		\begin{subfigure}{0.49\textwidth}
			\includegraphics[scale = 0.57]{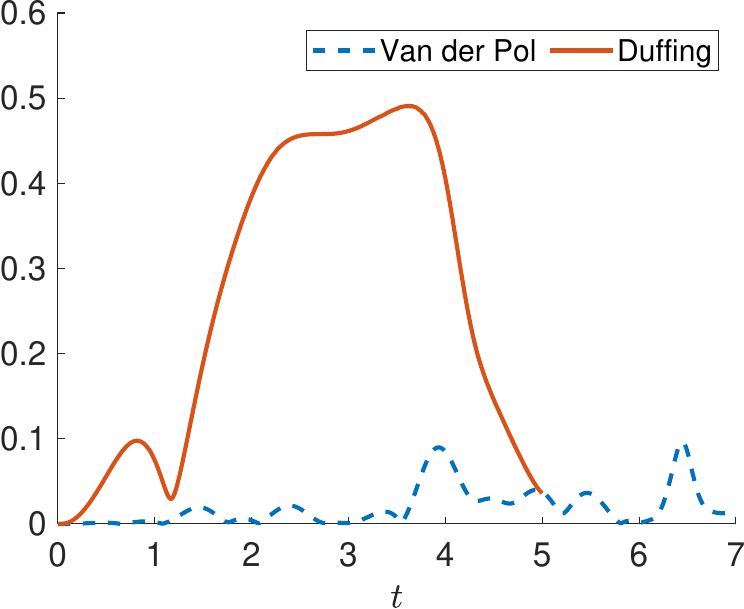}
			\caption{Relative differences of the observer gains}
			\label{fig: DiffGain}
		\end{subfigure}
		\caption{Comparing Van der Pol and Duffing}
		\label{fig: Compare VdP_Duff}
	\end{figure}

	A first attempt at investigating the cause of this behaviour is presented in \Cref{fig: Compare VdP_Duff}. There we compare the terms in which extended Kalman filter and Mortensen observer differ, c.f.~\Cref{sec: MorObs}. In \Cref{fig: DiffThirdDer} the ratio of the additional term in the Mortensen DRE and the right hand side of the extended Kalman filter DRE  
	\begin{equation*}
		\frac{\Pi(t) \nabla_{\xi^3}^3 \mathcal{V}(t,\widehat{x}_{\mathrm{M}}(t)) f(\widehat{x}_{\mathrm{M}}(t)) \Pi(t)}
		{\mathrm{D}f(\widehat{x}_{\mathrm{M}}(t)) \Pi(t) + \Pi(t) \mathrm{D}f(\widehat{x}_{\mathrm{M}}(t))^\top 
			- \alpha \Pi(t) C^\top C \Pi(t) + F F^\top}
	\end{equation*}
	is plotted over time. In \Cref{fig: DiffGain} we present the relative difference of the observer gains
	\begin{equation*}
		\frac{\Vert \Pi(t) C^\top - \Sigma(t) C^\top \Vert }{\Vert \Pi(t) C^\top \Vert}.
	\end{equation*} 
	
	Clearly  the quantifiers for the difference between the extended Kalman filter and the Mortensen observer are significant for the Duffing oscillator. They   give a first explanation for the noticeable difference in the state reconstruction based on these two methods. In the same figure we also present these quantifiers for the Van der Pol oscillator. It turns out that they are considerably smaller.   These observations certainly deserve further research.
	
	%figures observer
	
	%Figures third derivatives and gains

	%%%%%%%%%%%%%%%%%%%%%%%%%%%%%%%%%%%%%%%%%%%%%%%%%%%%%%%%
	\subsection{Test 4: Cucker-Smale model}\label{subs: test 4}~\\
	%%%%%%%%%%%%%%%%%%%%%%%%%%%%%%%%%%%%%%%%%%%%%%%%%%%%%%%%
	
	Finally we consider the Cucker-Smale model for agent based optimal consensus control. It should be noted that the original intent of modeling consensus behaviour is not considered in our discussion. We are interested in this model merely for its non-linear dynamics and the fact that the state space dimension can be easily adjusted by varying the number of agents. In order to apply our strategy to a system of medium sized state space dimension we consider the uncontrolled system with a disturbance. We discuss a system with $N_\mathrm{a}$ agents with states $(z_i,q_i) \in \mathbb{R}^2 \times \mathbb{R}^2$, for $i = 1,...,N_\mathrm{a}$. Here $z_i$ and $q_i$ correspond to the position and the velocity of the $i$-th agent moving in the plane. The dynamics are characterized by the equations
	\begin{equation*}
		\begin{aligned}
			\tfrac{\mathrm{d}}{\mathrm{d}t} z_i &= q_i \\
			\tfrac{\mathrm{d}}{\mathrm{d}t} q_i &= \frac{1}{N_\mathrm{a}} \sum_{j=1}^{N_\mathrm{a}} \frac{q_j - q_i}{1 + \Vert z_j - z_i \Vert^2} + v_i \\
			y_i &= z_i + \mu_i,
		\end{aligned}
	\end{equation*}
	where $v_i$ is the disturbance in the dynamics and $\mu_i$ represents the disturbance in the measurement. Analogous to the examples discussed above only the velocities are affected by system disturbances and the measurement consists of only the positions. The initial position and velocity of the $i$-th agent are set to 
	\begin{equation*}
		z_i(0) = q_i(0) = \frac{1}{2} 
		\begin{bmatrix}
			\cos (\frac{2 i \pi}{N_\mathrm{a}}) \\
			\sin (\frac{2 i \pi}{N_\mathrm{a}})
		\end{bmatrix}.
	\end{equation*}
	This choice places the agents on a circle around the origin and equips them with an initial velocity pointing straight away from the origin. For the construction of the measurement $y$ we set the disturbance in the velocity of the $i$-th agent to 
	\begin{equation*}
		v_i(t) = \frac{0.3}{\sqrt{N_\mathrm{a}}}
		\begin{bmatrix}
			0 & 1\\ -1 & 0
		\end{bmatrix}
		q_i(0).
	\end{equation*}
	For the error in the measurement we set 
	\begin{equation*}
		\mu_i(t) = \sin (2 \pi \, t) \frac{0.8}{\sqrt{N_\mathrm{a}}}
		\begin{bmatrix}
			0 & 1\\ -1 & 0
		\end{bmatrix}
		q_i(0).
	\end{equation*}
	We note that the direction of the disturbances is chosen to be orthogonal to the respective initial velocities. The specific choices for the initial condition and disturbances is not required to ensure satisfying results for our numerical scheme. Choosing random initial conditions from an appropriate domain and combining them with more general disturbances leads to approximations with the same order of accuracy. This particular setting was chosen because it leads to a visual representation with comparatively little overlap in the resulting trajectories.
	
\begin{figure}
		\centering
		\begin{subfigure}{0.45\textwidth}
			\includegraphics[scale = 0.5]{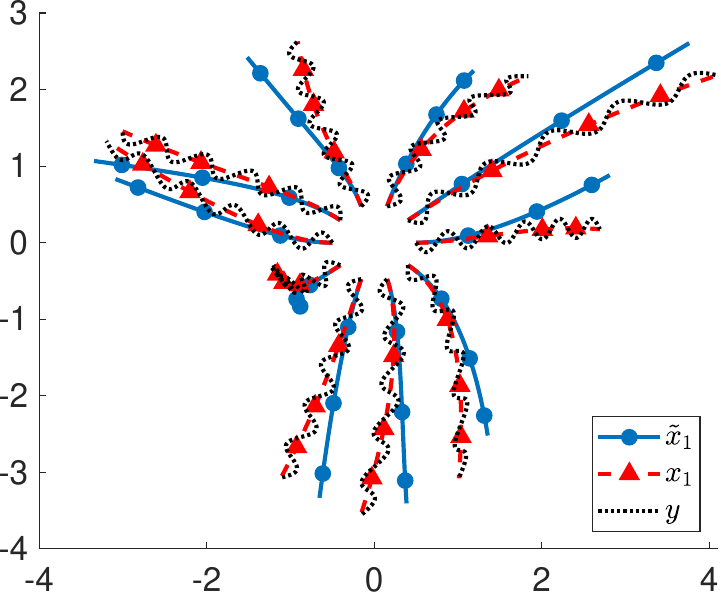}
			\caption{Positions of the model $\tilde{x}$, the disturbed trajectory $x$ and the measured positions $y$.}
			\label{fig: CuckerSmale_model_pos}
		\end{subfigure}
		\hfill
		\begin{subfigure}{0.45\textwidth}
			\includegraphics[scale = 0.5]{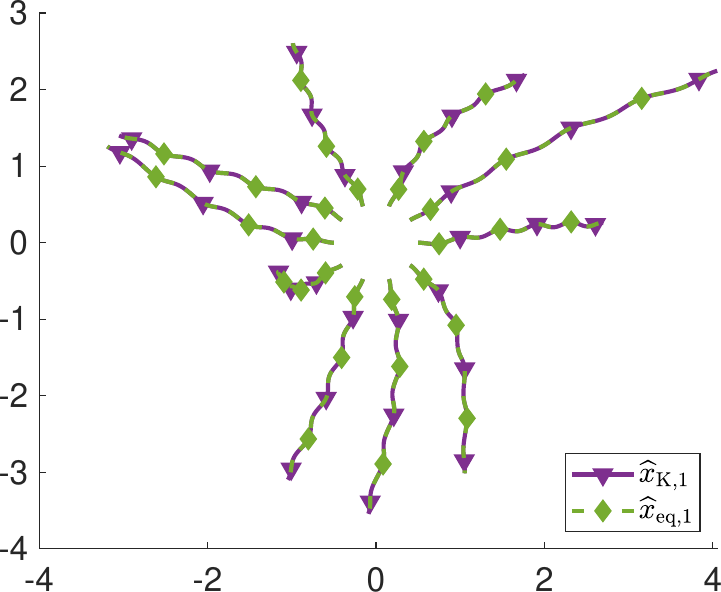}
			\caption{Positions of the extended Kalman filter $\widehat{x}_\mathrm{K}$ and the approximated Mortensen observer $\widehat{x}_\mathrm{eq}$.}
			\label{fig: CuckerSmale_obs_pos}
		\end{subfigure}
		\begin{subfigure}{0.45\textwidth}
			\includegraphics[scale = 0.5]{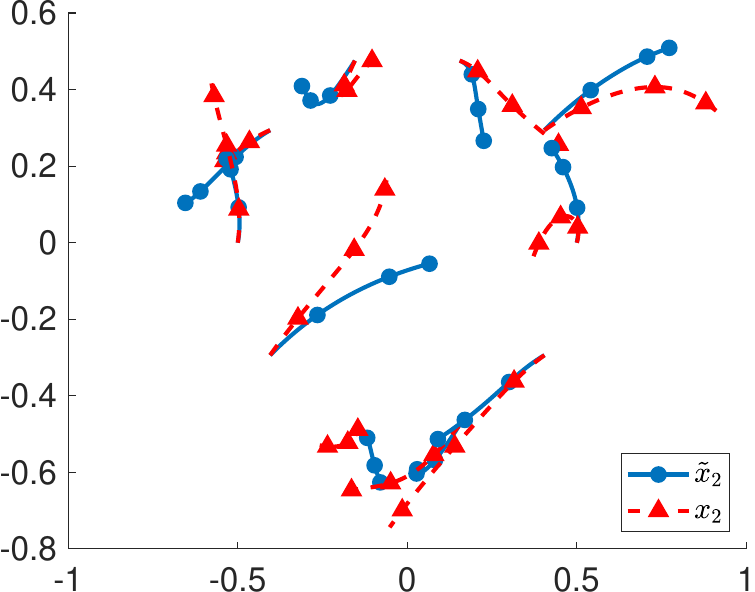}
			\caption{Velocities of the model $\tilde{x}$ and the disturbed trajectory $x$.}
			\label{fig: CuckerSmale_model_vel}
		\end{subfigure}
		\hfill
		\begin{subfigure}{0.45\textwidth}
			\includegraphics[scale = 0.5]{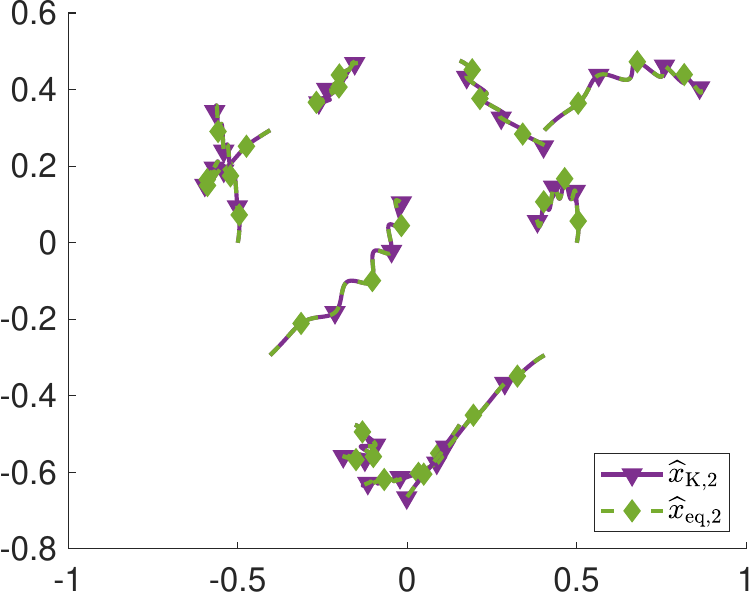}
			\caption{Velocities of the extended Kalman filter $\widehat{x}_\mathrm{K}$ and the approximated Mortensen observer $\widehat{x}_\mathrm{eq}$.}
			\label{fig: CuckerSmale_obs_vel}
		\end{subfigure}
		\caption{Cucker Smale equation}
		\label{fig: CuckerSmale}
	\end{figure}

	Again we focus only on the state reconstruction via the integration of the observer equation. For our computation we set the number of agents to $N_\mathrm{a} = 10$ resulting in a state space of dimension $n = 40$ and the time horizon is set to $[0,5]$. The domains for the spatial sampling are characterized by the side lengths $r_{k,i}$, where for $k=1,...,N_\mathrm{Time}$ and $i = 1,...,N_\mathrm{a}$ we set
	\begin{equation*}
		r_{k,2i} = r_{k,2i-1} = r_{k,2i+\tfrac{n}{2}} = r_{k,2i - 1 + \tfrac{n}{2}} 
		= \max \left\{ 0.1, \left\Vert 
		\begin{bmatrix}
			\widehat{x}_\mathrm{K}(t_k)_{2i-1} \\
			\widehat{x}_\mathrm{K}(t_k)_{2i}
		\end{bmatrix}
		\right\Vert \right\}.
	\end{equation*}
	We use $20$ time and $10$ space samples, respectively and the polynomial basis is constructed based on setting $d_\mathrm{Time} = 10$ and $s = 4$. The LS-problem considers the sampled values of $\mathcal{V}$ with a weight of $10^{-3}$ while the Hessians are entering with the weight $1$. The computations for this example took roughly 160 minutes. The resulting approximation of the Mortensen observer is compared with the extended Kalman filter in \Cref{fig: CuckerSmale_obs_pos} and \Cref{fig: CuckerSmale_obs_vel} and we observe that they agree. 
	
	%%%%%%%%%%%%%%%%%%%%%%%%%%%%%%%%%%%%%%%%%%%%%%%%%%%%%%%%%%%%%%
	\section{Conclusion}
	%%%%%%%%%%%%%%%%%%%%%%%%%%%%%%%%%%%%%%%%%%%%%%%%%%%%%%%%%%%%%%
	
	Two schemes for the numerical realization of the Mortensen observer were proposed.
	The examples under consideration show that both are viable options. In particular the integration of the observer equation based on an approximation of the value function yields satisfying results for systems of small and moderate state space dimension. The experiments further suggest that it is beneficial to not only consider samples of the values of the value function but also take its derivatives into account. Additionally, one of the examples shows that the Mortensen observer can substantially differ from the extended Kalman filter. 
	 
	%%%%%%%%%%%%%%%%%%%%%%%%%%%%%%%%%%%%%%%%%%%%%%%%%%%%%%%%%%%%%%
	%END OF MAIN ARTICLE%
	%%%%%%%%%%%%%%%%%%%%%%%%%%%%%%%%%%%%%%%%%%%%%%%%%%%%%%%%%%%%%%
	
	\section*{Acknowledgement}
	We thank B. H\"oveler (TU Berlin) for many helpful comments and suggestions on the practical aspects of solving ordinary differential equations. T. Breiten and J. Schr\"oder gratefully acknowledge funding and support from the Deutsche Forschungsgemeinschaft via the project 504768428.
	
	%%%%%%%%%%%%%%%%%%%%%%%%%%%%%%%%%%%%%%%%%%%%%%%%%%%%%%%%%%%%%%
	\bibliographystyle{siam}
	\bibliography{references_HessAugMor} 
	%%%%%%%%%%%%%%%%%%%%%%%%%%%%%%%%%%%%%%%%%%%%%%%%%%%%%%%%%%%%%%
	
\end{document}